\documentclass[twoside,11pt]{article}

\usepackage{fancyhdr}
\usepackage[colorlinks,citecolor=blue,urlcolor=blue,linkcolor=blue,bookmarks=false]{hyperref}
\usepackage{url}
\usepackage{doi}
\usepackage{epsf}
\usepackage{epsfig}
\usepackage{graphics}
\graphicspath{ {./figures/} }
\usepackage{float}
\usepackage{afterpage}
\usepackage{fullpage}
\usepackage[T1]{fontenc}
\usepackage{color,etoolbox}
\usepackage{enumitem}
\usepackage{import}
\usepackage[font=small,skip=0pt,labelfont=bf]{caption}
\usepackage[export]{adjustbox}
\usepackage{mdframed}

\usepackage{inconsolata}

\usepackage{tabulary}
\usepackage{booktabs}

\usepackage[square,comma,authoryear,sort&compress]{natbib}

\usepackage{./Latex/macros_01_proof_envs}
\usepackage{./Latex/macros_02_operators_commands}
\usepackage{./Latex/macros_03_symbols}
\usepackage{./Latex/macros_04_propernames}
\usepackage{./Latex/macros_05_reviewer}
\usepackage{./Latex/smile}

\usepackage[capitalize,sort,compress]{cleveref}

\usepackage[breakable]{tcolorbox}

\newtcolorbox{customcolorbox}[3][breakable]{colframe = #2!30, colback  = #2!8,
                                            coltitle = #2!30!black, title    =
                                            {#3}, #1, }

\newtcolorbox{mycolorbox}[3][breakable]{colframe = #2!25, colback  = #2!10,
                                        coltitle = #2!20!black, title    = {#3},
                                        #1, }

\usepackage{csquotes}
\newenvironment{itquote}
{\begin{quote}\itshape} {\end{quote}}

\newcommand{\diam}{\operatorname{diam}}

\def\T{{ ^\mathrm{\scriptscriptstyle T} }} \def\v{{\varepsilon}}

\newcommand{\hellmet}[2]{d_{\mathsf{H}}{\left(#1 || #2\right)}}
\newcommand{\kldiva}[2]{d_{\mathsf{KL}}{\left(#1 || #2\right)}}
\newcommand{\chisqdiva}[2]{\chi^{2}{\left(#1 || #2\right)}}

\newcommand{\cFalphabetaB}{\cF_{B}^{[\alpha, \beta]}}
\newcommand{\cFalphabetasupp}[1]{\cF_{ #1 }^{[\alpha, \beta]}}
\newcommand{\cFzerobetaB}{\cF_{B}^{[0, \beta]}}
\newcommand{\totbddvzeta}{\mathsf{BV}_{\zeta}}
\newcommand{\lipschitzpsi}{\mathsf{Lip}_{\gamma, q}(\Psi)}
\newcommand{\convclassk}{\mathsf{Conv}_{k}}
\newcommand{\quadfuncclass}{\mathsf{Quad}_{\gamma}}

\newcommand{\mgloc}[2]{M^{\mathsf{glo}}_{{#1}}\left(#2\right)}
\newcommand{\mlocc}[3]{M^{\mathsf{loc}}_{{#1}}\left(#2, #3\right)}
\newcommand{\madlocc}[4]{M^{\mathsf{adloc}}_{{#1}}\left(#2, #3, #4\right)}

\newcommand{\combdatavec}{\bfX}

\newtoggle{showfigs} 
\toggletrue{showfigs} 

\newtoggle{cmuthesismode}
\togglefalse{cmuthesismode} 

\begin{document}

\begin{center}
    {\LARGE \textbf{Revisiting Le Cam's Equation: Exact Minimax Rates over
        Convex Density Classes}}

    \vspace*{.2in}
    \begin{tabular}[t]{c@{\extracolsep{3em}}c} Shamindra Shrotriya  & Matey
    Neykov \\
    Retail Intelligence & Department of Statistics and Data Science \\
    Walmart & Northwestern University \\
    Bentonville, AR 72716 & Evanston, IL 60208 \\
    \texttt{shamindra.shrotriya@walmart.com} & \texttt{mneykov@northwestern.edu}
    \end{tabular}

    \vspace*{.1in}

    \iftoggle{cmuthesismode}{%
\begin{mdframed}[backgroundcolor=black!10,hidealllines=true]
  \textbf{Abstract:} 
  We study the classical problem of deriving minimax rates for density
estimation over convex density classes. Building on the pioneering work of
\citet{lecam1973convergence,birge1983approximation,
birge1986estimatinghellingerfacts,wong1995probability,yang1999information},
we determine the exact (up to constants) minimax rate over any convex
density class. This work thus extends these known results by demonstrating
that the local metric entropy of the density class always captures the
minimax optimal rates under such settings. Our bounds provide a unifying
perspective across both parametric and nonparametric convex density classes,
under weaker assumptions on the richness of the density class than
previously considered. Our proposed `multistage sieve' MLE applies to any
such convex density class. We further demonstrate that this estimator is also 
adaptive to the true underlying density of interest. We apply our risk bounds 
to rederive known minimax rates including bounded total variation, and Lipschitz 
density classes. We further illustrate the utility of the result by deriving upper
bounds for less studied classes, \eg, convex mixture of densities.\\
  \newline
  \nit The work in this chapter was done jointly with Matey Neykov. 
  It is based on a preprint with the title
  \emph{``Revisiting Le Cam's Equation: Exact Minimax Rates over
        Convex Density Classes''}.  
\end{mdframed}
}{%
\begin{abstract}
    \vspace*{.1in}
    We study the classical problem of deriving minimax rates for density
estimation over convex density classes. Building on the pioneering work of
\citet{lecam1973convergence,birge1983approximation,
birge1986estimatinghellingerfacts,wong1995probability,yang1999information},
we determine the exact (up to constants) minimax rate over any convex
density class. This work thus extends these known results by demonstrating
that the local metric entropy of the density class always captures the
minimax optimal rates under such settings. Our bounds provide a unifying
perspective across both parametric and nonparametric convex density classes,
under weaker assumptions on the richness of the density class than
previously considered. Our proposed `multistage sieve' MLE applies to any
such convex density class. We further demonstrate that this estimator is also 
adaptive to the true underlying density of interest. We apply our risk bounds 
to rederive known minimax rates including bounded total variation, and Lipschitz 
density classes. We further illustrate the utility of the result by deriving upper
bounds for less studied classes, \eg, convex mixture of densities.
\end{abstract}
}
\end{center}

\section{Introduction}\label{sec:introduction-density-estimation-convex-class}

It is well known that (global) metric entropy often times determines the minimax
rates for density estimation. Specifically, the following equation sometimes
informally referred to as the `Le Cam equation', is used to heuristically
determine the minimax rate of convergence
\begin{align*}
    \log{\mgloc{\cF}{\varepsilon}} \asymp n \varepsilon^2,
\end{align*}
where $n$ is the sample size, $\log{\mgloc{\cF}{\varepsilon}}$ is the
\emph{global} metric entropy of the density set $\cF$ at a Hellinger distance
$\varepsilon$ (see Definition \ref{nrmk:cF-compact-packing-numbers}), and
$\varepsilon^2$ determines the order of the minimax rate. In this paper we
complement these known results, by establishing that \emph{local} metric entropy
\emph{always} determines the minimax rate for convex density classes, where the
densities are assumed to be (uniformly) bounded from above and below.

In detail, under the setting of density estimation just described, we suggest a
small revision to the Le Cam  equation: namely, change the global entropy to
local entropy, and the Hellinger metric to the $L_2$-metric. Furthermore, the
same result holds when the convex density class contains densities only
(uniformly) bounded from above, and a single density which is bounded from
below. Unlike previous known results, our result unites minimax density
estimation under both parametric and nonparametric convex density classes. A
further contribution is that our proposed `multistage sieve' maximum likelihood
estimator (MLE) achieves these bounds regardless of the density class (as long
as it is convex). In addition, we demonstrate that this multistage sieve
estimator is adaptive to the true density of interest. To the best of our
knowledge we are not aware of any other estimators with such a property, under
such general settings.

We will now formally describe the setting we consider. To that end, we first
define a general class of bounded densities, \ie, $\cFalphabetaB$. Later, we
will assume that the true density of interest belongs to a known convex subset
of this general ambient density class.
\begin{restatable}[Ambient density class
        $\cFalphabetaB$]{ndefn}{densityclassFalphabeta}\label{ndefn:density-class-F-alpha-beta}
        Given constants $0 < \alpha < \beta < \infty$, for some fixed dimension
        $p \in \NN$, and a common known (Borel measurable) compact support set
        $B \subseteq \RR^{p}$ (with positive measure), we then define the class
        of density functions, $\cFalphabetaB$, as follows:
    \begin{equation}\label{neqn:density-class-F-alpha-beta}
        \cFalphabetaB
        \defined
        \thesetb{f \colon B \to [\alpha, \beta]}{\int_{B}{f} \dlett{\mu} = 1, \textnormal{$f$ measurable}},
    \end{equation}
    where $\mu$ is the dominating finite measure on $B$. We always take $\mu$ to
    be a (normalized) probability measure on $B$.
\end{restatable}

Furthermore, we can endow $\cFalphabetaB$ with the $L_{2}$-metric. That is, for
any two densities $f, g \in \cFalphabetaB$, we denote the $L_{2}$-metric between
them to be
\begin{equation}\label{neqn:L2-metric-density-defn}
    \normb{f - g}_{2}
    \defined
    \parens{\int_{B} (f - g)^{2} \dlett{\mu}}^{\frac{1}{2}}.
\end{equation}

\begin{restatable}{nrmk}{densityclasscfalphabeta}\label{nrmk:density-class-F-alpha-beta}
    Qualitatively, we have that $\cFalphabetaB$ is the class of all densities
    that are uniformly $\alpha$-lower bounded and $\beta$-upper bounded, on a
    common compact support $B \subseteq \RR^{p}$. Furthermore,
    \Cref{ndefn:density-class-F-alpha-beta} implies that $\cFalphabetaB$ forms a
    convex set, and that the metric space  $(\cFalphabetaB, \normb{\cdot}_{2})$
    is complete, bounded, but may not be totally bounded\footnote{These
    fundamental (and additional) analytic properties of $\cFalphabetaB$ are
    formally justified in \Cref{app:prop-of-class-F-alpha-beta}.}.
\end{restatable}

In this paper we will focus on the scenario where it is known that the true
density $f \in \cF \subset \cFalphabetaB$, where $\cF$ is a known convex set.
The set $\cF$ represents our knowledge on the true density, before observing any
data. With these mathematical preliminaries, we formalize our core density
estimation problem of interest as follows.
\begin{tcolorbox}
    \begin{itquote}
        \textbf{Core problem:}
        Suppose that we observe $n$ observations $\combdatavec \defined (X_1,
            \ldots, X_n)^{\top} \distiid f$, for some (fixed but unknown) $f \in
            \cF$. Here $\cF \subset \cFalphabetaB$ is a convex set, which is
            known to the observer. Can we propose a universal estimator for $f$,
            and derive the exact (up to constants) squared $L_{2}$-minimax rate
            of estimation, in expectation?
    \end{itquote}
\end{tcolorbox}

For convenience, we can illustrate the generating process for a univariate
example of our density estimation problem of interest in
\Cref{fig:generating-process}. It will serve as a useful conceptual guide to
later help visualize our proposed estimator over such general convex class of
densities $\cF$.

\iftoggle{showfigs}{%
    \begin{figure}[!ht]
        \centering
        \includegraphics[scale=0.75]{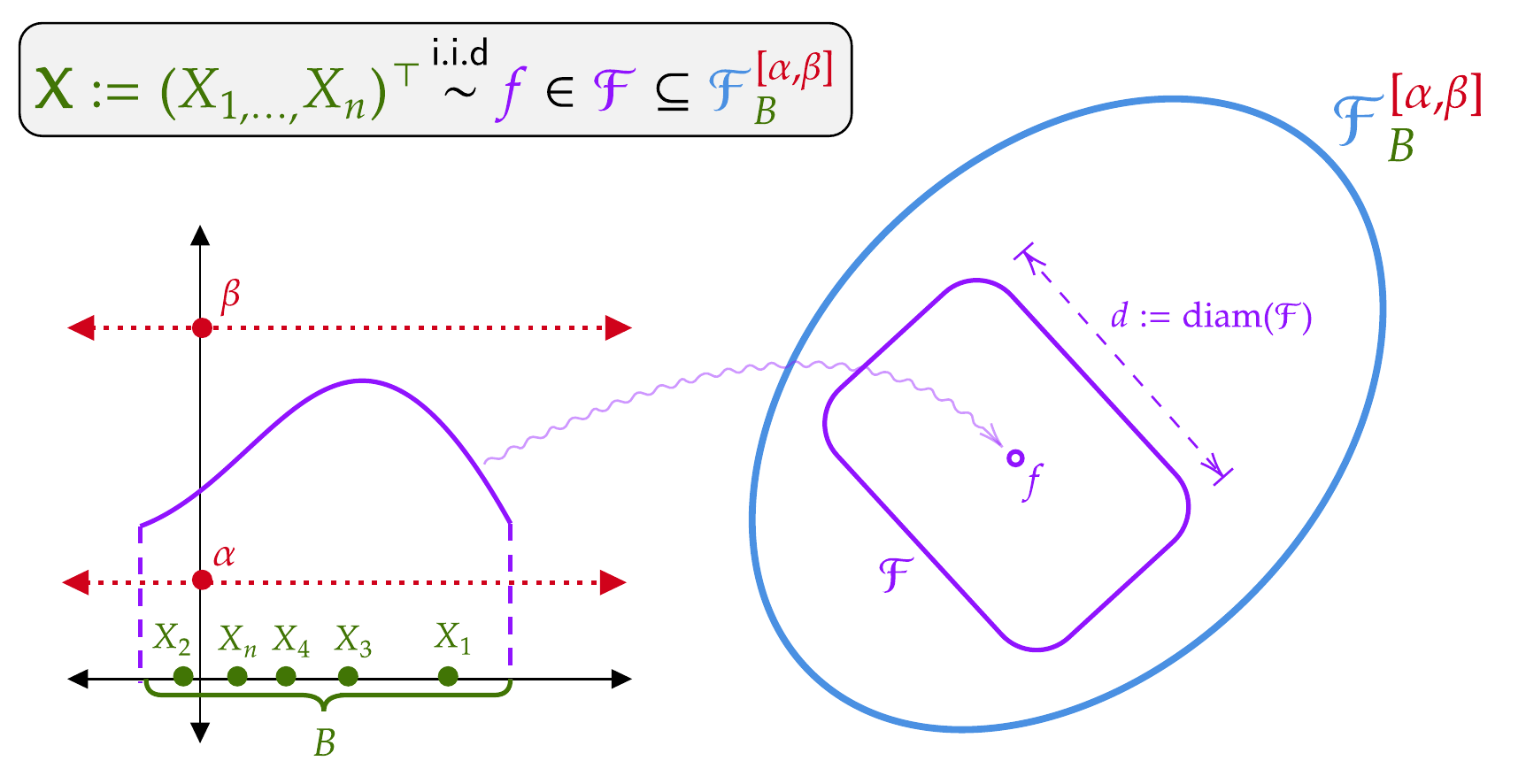}
        \caption{Illustrative generating process for a univariate density $f \in \cF \subset \cFalphabetaB$.}
        \label{fig:generating-process}
    \end{figure}
}{%
}

Now, without further ado, we will informally state our main result as a direct
affirmative answer to our core question of interest. Namely, there does exist a
likelihood-based estimator (one can think of it as a multistage sieve MLE), \ie,
$\nu^{*}(\combdatavec)$, which achieves the following rate of estimation error
\begin{equation}\label{neqn:intro-main-result-bound}
    \sup_{f \in \cF} \EE \normb{\nu^{*}(\combdatavec) - f}_{2}^{2} \lesssim \varepsilon^{*2} \wedge \diam_{2}(\cF)^2.
\end{equation}
Here $\varepsilon^* \defined \sup \thesetc{\varepsilon}{n \varepsilon^{2} \leq
        \log{\mlocc{\cF}{\varepsilon}{c}}}$, with
        $\log{\mlocc{\cF}{\varepsilon}{c}}$ being the $L_{2}$-\emph{local}
        metric entropy of $\cF$ (see Definition \ref{ndefn:local-entropy}). The
        quantity $\diam_{2}(\cF)$, refers to the $L_{2}$-diameter of $\cF$,
        which is finite by the boundedness of $\cFalphabetaB$ in our setting. In
        addition, the rate above is minimax optimal, as there is a matching (up
        to constants) lower bound.

\begin{restatable}{nrmk}{relaxingdensityclassFalphabeta}\label{nrmk:relaxing-density-class-F-alpha-beta}
    We will later see that we can largely relax the $\alpha$-lower boundedness
    condition on $\cF$. That is, the results we are about to derive can be
    readily generalized to convex subsets $\cF \subset \cFzerobetaB$. This is so
    as long as the class $\cF$ contains a \emph{single} density which is bounded
    away from $0$.
\end{restatable}

Next, we turn our attention to reviewing some relevant literature.

\subsection{Relevant Literature}\label{sec:relevant-literature}

\vspace{4pt}
\nit \underline{\textbf{Classical work}}

\vspace{4pt}
As noted, density estimation is a classical statistical estimation problem with
a rich history. Lively accounts of the key references, particularly for
nonparametric density estimation as relevant to our setting, are already covered
in \citet[Section~1]{yang1999information} and
\citet[Section~6.1]{bilodeau21density}. We similarly begin with a brief
panoramic overview of these references in regard to minimax risk bounds for
density estimation, before comparing and contrasting the results from the most
relevant references to our work.

In terms of minimax lower bounds on density estimation,
\citet{boyd1978lowerbddsnonparamest} prove a fundamental $n^{-1}$ rate in the
mean integrated $p^{\textnormal{th}}$ power error (with $p \geq 1$), for
\emph{any} arbitrary density estimator. Such generalized lower bounds on density
estimation were also further studied in \citet{devroye1983arbslowratesconv}. In
the case of density estimation over classes with more assumed structure (\eg,
smoothness, or regularity assumptions) minimax lower bounds have been developed
based on hypothesis testing approaches coupled with information-theoretic
techniques. We now provide brief highlights of such key works in this direction.

In \citet{bretagnolle1979estrisqueminimax}, the authors derive sharp lower
bounds for Sobolev smooth densities in $\RR^{d}$ ($d \in \NN$), with risk
measured with respect to a power of the $L_{q}$-metric ($q \geq  1$). In
\citet{birge1986estimatinghellingerfacts}, sharp risk bounds for more general
classes of such smooth families were provided using metric entropy based
methods, with an emphasis on the Hellinger loss. The work of
\citet{efroimovich1982estsquareintdensityrv} provided precise (asymptotic)
analysis for an ellipsoidal class of densities in the $L_{2}$-metric. Across a
wide-ranging series of related and collaborative efforts
\citet{hasminski1978lowerboundunifmetric, ibragimov1977estimationparwhitenoise,
ibragimov1978capacitycommsmoothsig,ibragimov1980estdistdensity} used Fano's
inequality type arguments to establish lower bounds over a variety of density
estimation settings. These range from deriving lower bounds on nonparametric
density estimation in the uniform metric, to minimax risk bounds for the
Gaussian white noise model, for example. The authors also develop metric entropy
based techniques in \cite{hasminski1990denskolmogorovappr} to derive minimax
lower bounds for a wide variety of density classes defined on $\RR^{d}$ ($d  \in
\NN$), in $L_{q}$-loss ($q \geq 1$). Numerous applications of optimal lower
bounds using both Assouad's and Fano's lemma arguments for densities on a
compact support, are demonstrated in
\citep[Section~29.3]{yu1997assouadfanolecamminimaxlb}. Later
\citet{yang1999information} demonstrated that \emph{global} metric entropy
bounds capture minimax risk for sufficiently rich density classes over a common
compact support. Classical reference texts on minimax lower bound techniques
with an emphasis on nonparametric density estimation include
\citet{devroye1987coursedensityest,devroye1985nonparamdensityestl1view,
lecam1986asympmethodsdectheory}. More modern such references include
\citet{tsybakov2009introduction} and \citet[Chapter~15]{wainwright2019high}. The
latter in particular, also incorporates metric entropy based lower bound
techniques.

In addition, there is a large body of work in deriving upper bounds for specific
density estimators using metric entropy methods. This includes
\citet{yatracos1985ratesconvmindistent, barron1991mincompdensityest}, who employ
the minimum distance principle to derive density estimators and their metric
entropy-based upper bounds in the Hellinger and $L_{1}$-metric, respectively. In
a similar spirit to \citet{birge1983approximation,
birge1986estimatinghellingerfacts}, \citet{vandegeer1993hellingerconsnpmle} is
also concerned with density estimation using Hellinger loss. However, its focus
is to use techniques from empirical process theory in order to specifically
establish the Hellinger consistency of the nonparametric MLE, over convex
density classes. Upper bounds for density estimation based on the `sieve' MLE
technique is studied in \citet{wong1995probability}. Recall that a `sieve'
estimator effectively estimates the parameter of interest via an optimization
procedure (\eg, maximum likelihood) over a constrained subset of the parameter
space \citep[Chapter~8]{grenander1981abstractinference}. In
\citet{birge1993ratesconvmincontrast} the authors study `minimum contrast
estimators' (MCEs), which include the MLE, least squares estimators (LSEs) \etc,
and apply them to density estimation. This is further developed in
\cite{birge1994mincontrastsieve} where they analyze convergence of MCEs using
sieve-based approaches.

\vspace{4pt}
\nit \underline{\textbf{Comparison to our work}}

\vspace{4pt}
By stating our main result early in the introduction, we now turn to contrasting
it with the most relevant results in the literature. These include both the
aforementioned classical references, and more recent work on convex density
estimation, which have most directly inspired our efforts in this work.

First we would like to comment on the closely related landmark papers
\citep{lecam1973convergence,
birge1983approximation,birge1986estimatinghellingerfacts}. These works consider
very abstract settings and show upper bounds based on Hellinger ball testing.
Although widely believed that they do, whether these results lead to bounds that
are minimax optimal is unclear. Moreover, their estimator is quite involved and
non-constructive. In contrast, in this paper we offer a simple to state,
\emph{constructive} multistage sieve MLE type of estimator, which is provably
minimax optimal over any convex density class $\cF$. A crucial difference is
that we metrize the space $\cF$ with the $L_{2}$-metric as we mentioned above.
Even though in our instance the two distances are equivalent, in contrast to the
Hellinger distance, the $\varepsilon$-\emph{local} metric entropy of the convex
density class in the $L_{2}$-metric can be shown to be monotonic in
$\varepsilon$. This key observation enables us to match the upper and lower
bounds exactly.

Next, we will compare our work with the celebrated paper of
\cite{yang1999information}, who inspect a very similar problem.
\cite{yang1999information} demonstrate a lower and upper bound which need not
match in general but do match under certain sufficient conditions. Notably their
bounds involve only quantities depending on the global entropies of the set
$\cF$ (which is also assumed to be convex for some results of
\cite{yang1999information}). This is convenient as often times global metric
entropy is easier to work with compared to local metric entropy, however under
\cite{yang1999information}'s sufficient condition it can be seen that the two
notions are equivalent. Hence, our work can be thought of as removing the
sufficient condition requirement from \cite{yang1999information} and also
unifying parametric and nonparametric density estimation problems (over convex
classes) for which one typically needs to use different tools to obtain the
accurate rates. Finally we would like to mention \cite{wong1995probability}. In
that paper the authors propose a sieve MLE estimator and demonstrate that it is
\emph{nearly} minimax optimal under certain conditions. Our estimator is not the
same as the one considered by \cite{wong1995probability}, and we can provably
match the minimax rate over whatever be the convex set $\cF$. A notable
difference is that \cite{wong1995probability} work with the Hellinger metric and
KL divergence, which although equivalent to $L_{2}$-metric in our problem, are
actually less practical in terms of matching the bounds exactly as we explained
above. We will now turn our attention to reviewing some further relevant
literature.

\vspace{4pt}
\nit \underline{\textbf{Recent work}}

\vspace{4pt}

Our estimator and proof techniques thereof, are inspired by the recent work of
\cite{neykov2022minimax} on the Gaussian sequence model. We would like to stress
on the fact that the sequence model is a very distinct problem from density
estimation. In particular, our underlying metric space of interest is $(\cF,
\normb{\cdot}_{2})$, as compared to\footnote{Note that $\normb{\cdot}_{2}$  here
is the Euclidean metric on $\RR^{n}$.} $(\RR^{n}, \normb{\cdot}_{2})$ for the
sequence model. Both of these metric spaces differ \emph{vastly} from each other
in their underlying geometric structure. Furthermore, unlike our setting, the
sequence model contains additional Gaussian information on the underlying
generating process, which can be directly exploited for estimation purposes. As
such, given that \cite{neykov2022minimax} provides a guiding template for our
analysis, some resulting structural similarities to their work are to be
expected. However, all corresponding proofs, and estimators thereof, have to be
non-trivially adapted to our nonparametric density estimation setting. A notable
example of such required modifications, is that our estimator presented in this
paper does not use proximity in Euclidean norm, but is a likelihood-based
estimator.

We additionally note that density estimation in both abstract and more concrete
settings, continues to be  an active area of research. It is not feasible to
detail such a large and growing body of references. However, we provide a
selective overview of some interesting recent directions in density estimation,
to simply indicate the diversity of the research efforts thereof. For example,
\citet{cleanthous2020kernwaveletdensityest, baldi2009adaptdensityestspherical}
study convergence properties of density estimators using wavelet-based methods.
The papers
\citet{goldenshluger2014onadaptivedensityestrd,efromovich2008adaptiveestoracle,
rigollet2006adaptiveestblockstein,rigollet2007linearconvaggdensityest,
samarov2007aggdensityest,birge2014modselectdensityest} study adaptive minimax
density estimation on $\RR^{d}$ ($d \geq 1$) under $L_{p}$-loss ($p \geq 1$).
Here, `adaptive' refers to the fact that the density class is defined by an
unknown tuning hyperparameter, which must be explicitly accounted for during the
estimation process. Recently \citet{wang2022minimaxdensityoptimaltransport} used
techniques from optimal transport to study the convergence properties of various
nonparametric density estimators. Interestingly, \citet{bilodeau21density}
applied empirical (metric) entropy methods to establish minimax optimal rates in
the adjacent setting of \emph{conditional} density estimation. Although these
works do not directly study our core problem of interest, we note that they
represent new and important perspectives on classical minimax density
estimation, and related problems.


\subsection{Notation}\label{sec:notation-density-estimation-convex-class}

We outline some commonly used notation here. We use $a \vee b$ and $a \wedge b$
for the $\max$ and $\min$ of two numbers $\theset{a, b}$, respectively.
Throughout the paper $\normb{\cdot}_{2}$ denotes the $L_{2}$-metric in $\cF$.
Constants may change values from line to line. For an integer $m\in \NN$, we use
the shorthand $[m] \defined \{1, \ldots, m\}$. We use $B_2(\theta,r)$ to denote
a closed $L_{2}$-ball centered at the point $\theta$ with positive radius $r$.
We use $\lesssim$ and $\gtrsim$ to mean $\leq$ and $\geq$ up to absolute
(positive) constant factors, and for two sequences $a_n$ and $b_n$ we write $a_n
\asymp b_n$ if both $a_n \lesssim b_n$ and $a_n \gtrsim b_n$ hold. Throughout
the paper we use $\log$ to denote the natural logarithm, or we specify the base
explicitly otherwise. Our use of $\theset{\alpha, \beta}$ is \emph{only} used to
refer to the constants in \Cref{ndefn:density-class-F-alpha-beta}, of
$\cFalphabetaB$ (and thus $\cF$). We will introduce additional section-specific
notation as needed.

\subsection{Organization}\label{sec:organization-paper} The rest of this paper
is organized as follows. In \Cref{sec:minimax-upper-and-lower-bounds} we prove
risk bounds for our underlying setting. We first establish the key topological
equivalence between the $L_{2}$-metric and the Kullback-Leibler divergence in
$\cFalphabetaB$. We then proceed to derive minimax lower bounds for our setting
in \Cref{sec:lower-bound}, introducing additional relevant mathematical
background as needed, \eg, local metric entropy. In \Cref{sec:upper-bound} we
define our likelihood-based estimator, and provide intuition behind its
construction. We then derive its (matching) minimax risk upper bound. We then
demonstrate in \Cref{sec:adaptivity} that our estimator is adaptive to the true
density $f$. In \Cref{sec:examples}, we apply our results to specific examples
of commonly used convex density classes. We then conclude in
\Cref{sec:discussion} by summarizing our results, and discuss some future
research directions.

\section{Minimax Lower and Upper
  Bounds}\label{sec:minimax-upper-and-lower-bounds}

Before establishing our main results, we establish a key technical lemma which
drives much of the geometric arguments in our analysis to follow. Note that for
any two densities $f, g \in \cFalphabetaB$, the $\mathsf{KL}$-divergence between
them is defined to be
\begin{equation}\label{neqn:kldiv-density-defn}
    \kldiva{f}{g}
    \defined
    \int_{B} f \log\parens{\frac{f}{g}} \dlett{\mu}
    \defines
    \EE_{f}\log\parens{\frac{f(X)}{g(X)}},
\end{equation}
where $X \sim f$ in \eqref{neqn:kldiv-density-defn}.
\begin{restatable}{nrmk}{kldivdensitywelldefined}\label{nrmk:kldiv-density-well-defined}
    We observe that $\kldiva{f}{g}$ is well-defined in
    \eqref{neqn:kldiv-density-defn} for our setting, since $\inf_{x \in B} g(x)
    \geq \alpha > 0$, by \Cref{ndefn:density-class-F-alpha-beta}. We further
    emphasize that $\mathsf{KL}$-divergence is not valid metric in general,
    since it is not symmetric in its arguments.
\end{restatable}
The crucial fact in the risk bounds we will soon derive, is the `topological
equivalence' of the $L_{2}$-metric and $\mathsf{KL}$-divergence, on the density
class $\cFalphabetaB$. Since it is hard to find a concrete reference for this
folklore fact, we formalize this equivalence for our setting in
\Cref{nlem:equiv-kl-euc-metric}.
\begin{restatable}[$\mathsf{KL}$-$L_{2}$ equivalence on
        $\cFalphabetaB$]{nlem}{equivkleucmetric}\label{nlem:equiv-kl-euc-metric}
        For each pair of densities $f, g \in \cFalphabetaB$, the following
        relationship holds:
    \begin{equation}\label{neqn:nlem:equiv-kl-euc-metric-01}
        c(\alpha,\beta) \normb{f-g}_{2}^{2}
        \leq \kldiva{f}{g}
        \leq (1 / \alpha) \normb{f - g}_{2}^{2},
    \end{equation}
    where we denote $c(\alpha,\beta) \defined \frac{h(\beta / \alpha)}{\beta} >
        0$. Here $h : (0, \infty) \to \RR$ is defined to be
    \begin{equation}\label{neqn:nlem:equiv-kl-euc-metric-02}
        h(\gamma)
        \defined
        \begin{cases}
            \frac{\gamma - 1 - \log{\gamma}}{(\gamma - 1)^{2}} & \text{if $\gamma \in (0, \infty) \setminus \theset{1}$} \\
            \frac{1}{2} = \lim_{x \to 1} \frac{x - 1 - \log{x}}{(x - 1)^{2}}
                                                               & \text{if $\gamma = 1$},
        \end{cases}
    \end{equation}
    and is positive over its entire support. It is also easily seen that on
    $\cFalphabetaB$, $d_{\mathsf{KL}}$ (and hence the $L_{2}$-metric) is also
    equivalent to the Hellinger metric. Furthermore, these properties are also
    inherited by $\cF \subset \cFalphabetaB$, which is our density class of
    interest.
\end{restatable}
\begin{restatable}{nrmk}{klemelaequivkleucmetric}\label{nrmk:equiv-kl-euc-metric}
    We note that both the upper and lower bounds in
    \eqref{neqn:nlem:equiv-kl-euc-metric-01} are stated without proof and
    without tracking constants in \citet[Lemma~11.6]{klemela2009smoothing}. We
    formally prove this claim in \Cref{app:mathematical-preliminaries}.
    Importantly, the validity of \eqref{neqn:nlem:equiv-kl-euc-metric-01} relies
    on the assumption of the boundedness of the densities, which holds in our
    setting.
\end{restatable}

\subsection{Minimax Lower Bound}\label{sec:lower-bound} We will first establish
a lower bound. For completeness, we need to introduce some additional relevant
background and notation. We start by stating Fano's inequality for our convex
density class, $\cF$ \citep[see][Lemma~2.10]{tsybakov2009introduction}.
\begin{restatable}[Fano's inequality for
        $\cF$]{nlem}{fanoinequality}\label{nlem:fano-inequality} Let
        $\theset{f^1, \ldots, f^m} \subset \cF$ be a collection of
        $\varepsilon$-separated densities (\ie $\normb{f^{i} - f^{j}}_{2} >
        \varepsilon$ for $i \neq j$), in the $L_{2}$-metric. Suppose $J$ is
        uniformly distributed over the index set $[m]$, and $(X_{i} | J = j)
        \distiid f^j$ for each $i \in [n]$. Then
    \begin{align*}
        \inf_{\widehat{\nu}} \sup_{f} \EE \|\widehat{\nu}(\combdatavec) - f\|_{2}^{2}
        \geq
        \frac{\varepsilon^2}{4}\bigg(1 - \frac{n I(X_1; J) + \log 2}{\log m}\bigg).
    \end{align*}
\end{restatable}
In the above $I(X_1;J) \defined \frac{1}{m}\sum_{j = 1}^{m} \kldiva{f^{j}}{\bar
        f}$, where $\bar f = \frac{1}{m}\sum_{j = 1}^m f^j$ is the mutual
        information between $X_1$ and the randomly sampled index $J$. Further,
        the infimum is taken over all measurable functions of the data. Next, we
        define the important notion of a packing set for $\cF$ \citep[see
        Section~5.2][\eg, for more details]{wainwright2019high}.
\begin{restatable}[Packing sets and packing numbers of $\cF$ in the
        $L_{2}$-metric]{ndefn}{packingsets}\label{ndefn:packing-sets} Given any
        $\varepsilon > 0$, an $\varepsilon$-packing set of $\cF$ in the
        $L_{2}$-metric, is a set $\theset{f^1, \ldots, f^m} \subset \cF$ of
        $\varepsilon$-separated densities (\ie, $\normb{f^{i} - f^{j}}_{2} >
        \varepsilon$ for $i \neq j$) in the $L_{2}$-metric. The corresponding
        $\varepsilon$-packing number, denoted by $M(\varepsilon, \cF)$, is the
        cardinality of the largest (maximal) $\varepsilon$-packing of $\cF$. We
        refer to $\log{\mgloc{\cF}{\varepsilon}} \defined \log{M(\varepsilon,
        \cF)}$ as the \emph{global metric entropy} of $\cF$.
\end{restatable}
\begin{restatable}{nrmk}{cFcompactpackingnumbers}\label{nrmk:cF-compact-packing-numbers}
    Note that we are not assuming here that $\cF$ is totally bounded, hence some
    (or perhaps all) of the packing numbers may be infinite; this however does
    not cause a problem in what follows. Henceforth, all packing sets (or
    packing numbers) of $\cF$, will be assumed to be with reference to the
    $L_{2}$-metric, unless stated otherwise. We will use the standard fact that
    a $\varepsilon$-maximal packing of $\cF$, is also a $\varepsilon$-covering
    set of $\cF$.
\end{restatable}

We will now define the notion of \emph{local metric entropy}, which will play a
key role in the development of our risk bounds.
\begin{restatable}[Local metric entropy of
        $\cF$]{ndefn}{localentropy}\label{ndefn:local-entropy} Let $c > 0$ be
        fixed, and $\theta \in \cF$ be an arbitrary point. Consider the
        set\footnote{Observe that this set may also fail to be totally bounded,
        since while the ball $B_2(\theta,\varepsilon)$ is a bounded set, it is
        not totally bounded.} $\cF \cap B_{2}(\theta, \varepsilon)$. Let
        $M(\varepsilon/c, \cF \cap B_{2}(\theta, \varepsilon))$ denote the
        $\varepsilon/c$-packing number of $\cF \cap B_{2}(\theta, \varepsilon)$,
        in the $L_{2}$-metric. Let
    \begin{align*}
        \mlocc{\cF}{\varepsilon}{c}
        \defined
        \sup_{\theta \in \cF} M(\varepsilon/c, \cF \cap B_{2}(\theta, \varepsilon))
        \defines
        \sup_{\theta \in \cF} \mgloc{\cF \cap B_{2}(\theta, \varepsilon)}{\varepsilon/c}.
    \end{align*}
    We refer to $\log{\mlocc{\cF}{\varepsilon}{c}}$ as the \emph{local metric
        entropy} of $\cF$.
\end{restatable}

We show the following minimax lower bound for our convex density estimation
setting over $\cF$. It is a direct consequence of Fano's inequality per
\Cref{nlem:fano-inequality}.

\begin{restatable}[Minimax lower
        bound]{nlem}{minimaxlowerbound}\label{nlem:minimax-lower-bound} Let $c >
        0$ be fixed, and independent of the data samples $\combdatavec$. Then
        the minimax rate satisfies
    \begin{align*}
        \inf_{\widehat{\nu}} \sup_{f \in \cF} \EE_f \normb{\widehat{\nu}(\combdatavec) - f}_{2}^{2} \geq \frac{\varepsilon^2}{8c^2},
    \end{align*}
    if $\varepsilon$ satisfies $\log{\mlocc{\cF}{\varepsilon}{c}} > 2 n
        \varepsilon^{2} / \alpha + 2\log 2$.
\end{restatable}

\subsection{Upper Bound}\label{sec:upper-bound}

We now turn our attention to the upper bound. We note that our universal
estimator over $\cF$, will be a likelihood-based estimator for $f$. As such for
any two densities $g, g^{\prime} \in \cF$, we will routinely work with the
\emph{log-likelihood difference} for the $n$ observed samples $\combdatavec
\defined (X_{1}, \ldots, X_{n})^{\top} \distiid f \in \cF$. We will denote this
by
\begin{equation}\label{neqn:log-likelihood-g-g-prime}
    \psi(g, g^{\prime}, \combdatavec)
    \defined
    \log{\parens{\prod_{i = 1}^{n} \frac{g(X_{i})}{g^{\prime}(X_{i})}}}
    =
    \sum_{i = 1}^{n} \log\parens{\frac{g(X_{i})}{g^{\prime}(X_{i})}}
    =
    \sum_{i = 1}^{n} \log g(X_{i}) - \sum_{i = 1}^{n} \log g^{\prime}(X_{i}).
\end{equation}
\begin{restatable}{nrmk}{loglikelihoodwelldefined}\label{nrmk:log-likelihood-well-defined}
    We note that the log-likelihood difference $\psi(g, g^{\prime},
    \combdatavec)$ in \eqref{neqn:log-likelihood-g-g-prime}, is well-defined.
    This follows since for each $i \in [n]$, the individual random variables
    $\log g(X_{i})/g^{\prime}(X_{i})$ are well-defined (as $\alpha > 0$), and
    bounded. That is, $-\infty < \log \alpha/\beta \leq \log
    g(X_{i})/g^{\prime}(X_{i}) \leq \log \beta/\alpha < \infty$, for each $i \in
    [n]$.
\end{restatable}
We will use the log-likelihood difference to help us decide which of the two
densities is ``more'' correct, given the observed data samples $\combdatavec$.
Given this, we will first need a concentration result on the density
log-likelihood difference. We do this by establishing the following lemma.

\begin{restatable}[Log-likelihood difference concentration in
        $\cF$]{nlem}{loglikelihoodhoeffding}\label{nlem:log-likelihood-bernstein}
        Let $\delta > 0$ be arbitrary, and let $\combdatavec \defined (X_{1},
        \ldots, X_{n})^{\top} \distiid f \in \cF$, be the $n$ observed samples.
        Suppose we are trying to distinguish between two densities $g,
        g^{\prime} \in \cF$. Let $\psi(g, g^{\prime}, \combdatavec)$ denote
        their log-likelihood difference per
        \eqref{neqn:log-likelihood-g-g-prime}. We then have
    \begin{equation}\label{neqn:log-likelihood-bernstein-01}
        \sup_{\substack{g, g^{\prime} \colon \normb{g - g^{\prime}}_{2} \geq C\delta, \\
                \normb{g^{\prime}-f}_{2} \leq \delta}}\PP(\psi(g, g^{\prime}, \combdatavec) > 0)
        \leq
        \exp\parens{-n L(\alpha, \beta, C) \delta^{2}}
    \end{equation}
    where
    \begin{align}
        C
         & >
        1 + \sqrt{1 / (\alpha c(\alpha,\beta))}
        \label{neqn:log-likelihood-bernstein-02a} \\
        L(\alpha, \beta, C)
         & \defined
        \frac{\parens{\sqrt{c(\alpha,\beta)} (C-1)  -  \sqrt{1 / \alpha}}^{2}
        }{2\braces{ 2 K(\alpha, \beta) +\frac{2}{3} \log \beta / \alpha}},
        \label{neqn:log-likelihood-bernstein-02b}
    \end{align}
    with $K(\alpha, \beta) \defined \beta / (\alpha^{2} c(\alpha, \beta))$, and
    $c(\alpha, \beta)$ is as defined in \Cref{nlem:equiv-kl-euc-metric}. In the
    above $\PP$ is taken with respect to the true density function $f$, \ie,
    $\PP = \PP_f$.
\end{restatable}

From \Cref{nlem:log-likelihood-bernstein}, we derive a key concentration result
concerning a packing set in $\cF$, as summarized in
\Cref{nlem:critical-log-likelihood-concentration}. The relevance of such a
result will become clearer later, when we introduce our sieve-based MLE for $f$.
Our sieve estimator will be constructed using packing sets of $\cF$, thus
\Cref{nlem:critical-log-likelihood-concentration} will be an important tool to
enable us to handle the concentration properties of our estimator.

\begin{restatable}[Maximum likelihood concentration in
        $\cF$]{nlem}{criticalloglikelihoodconcentration}\label{nlem:critical-log-likelihood-concentration}
        Let $\delta > 0$ be arbitrary, and let $\combdatavec \defined (X_{1},
        \ldots, X_{n})^{\top} \distiid f \in \cF$, be the $n$ observed samples.
        Suppose further that we have a maximal $\delta$-packing set of
        $\cF^{\prime} \subset \cF$, \ie, $\theset{g_1, \ldots, g_m} \subset
        \cF^{\prime}$ such that $\normb{g_i- g_{j}}_{2} > \delta$ for all $i
        \neq j$, and it is known that $f \in \cF^{\prime}$. Now let $j^{*} \in
        [m]$, denote the index of a density whose likelihood is the largest. We
        then have
    \begin{equation*}
        \PP(\|g_{j^*}-f\|_2 > (C + 1) \delta) \leq
        m \exp\parens{-n L(\alpha, \beta, C) \delta^{2}},
    \end{equation*}
    where $C$ is assumed to satisfy \eqref{neqn:log-likelihood-bernstein-02a},
    and $L(\alpha, \beta, C)$ is defined as per
    \eqref{neqn:log-likelihood-bernstein-02b}.
\end{restatable}

Next we establish that the map $\varepsilon \mapsto \log
    \mlocc{\cF}{\varepsilon}{c}$ is non-increasing. This lemma is made possible
    by the fact that the set $\cF$ is convex by assumption, and that we are
    using the $L_{2}$-metric. This monotonicity property of the
    $\varepsilon$-local metric entropy in the $L_{2}$-metric is a
    \emph{critical} technical ingredient used in the proofs establishing our
    upper bound.
\begin{restatable}[Monotonicity of local metric
        entropy]{nlem}{localmetricentmonotone}\label{nlem:local-metric-ent-monotone}
        The map $\varepsilon \mapsto \log{\mlocc{\cF}{\varepsilon}{c}}$ is
        non-increasing.
\end{restatable}

We now turn our attention to describing our proposed likelihood-based estimator,
\ie, $\nu^{*}(\combdatavec)$, of $f \in \cF$. In the discussion that follows we
let $d \defined \diam_{\operatorname{2}}(\cF)$, which is finite by the
boundedness of $\cF$. The estimator is directly inspired by a recent
construction used in \citet{neykov2022minimax}, who applied it to the Gaussian
sequence model. However, there the underlying space used is $(\RR^{n},
\normb{\cdot}_{2})$, whereas in our case it is $(\cF, \normb{\cdot}_{2})$, which
has a \emph{vastly} different underlying geometric structure. Importantly since
we are performing density estimation, our proposed estimator uses a
fundamentally different \emph{log-likelihood}-based selection criteria, compared
to the \emph{projection}-based sequence model estimator in
\citet{neykov2022minimax}. Similar to \citet{neykov2022minimax}, although our
estimator can also be described constructively, it is not intended to be
practically computable. Our estimator will be shortly described as a four-stage
constructive procedure (\ie,
\ref{itm:likelihood-based-estimator-01}--\ref{itm:likelihood-based-estimator-04}).
Since the qualitative description of the construction may appear to be quite
involved, we provide a corresponding visual representation in
\Cref{fig:packing-set-tree-construction} as a helpful guide for the reader. We
emphasize however, that \Cref{fig:packing-set-tree-construction} is not drawn to
any precise scale.

\vspace{8pt}
\begin{mycolorbox}{blue}{\textbf{Construction of the multistage sieve MLE,
            $\nu^{*}(\combdatavec)$, of $f \in \cF$.}}
    \benum[wide=0pt, label={\color{blue} \textbf{Step \arabic*}}, align=left, start=1]
    \item \label{itm:likelihood-based-estimator-01} \textbf{Initialize inputs.}
    \vspace{1mm}
    \newline
    \nit
    Let $\combdatavec \defined (X_{1}, \ldots, X_{n})^{\top}$ denote our $n$
    observed \iid data samples. Fix some sufficiently large $c > 0$, and then
    define $C$ such that $c \defined 2(C + 1)$. Importantly, the constant $c$
    should be set \emph{without} looking at the data samples, \ie,
    \emph{independently} of $\combdatavec$.
    \item \label{itm:likelihood-based-estimator-02} \textbf{Construct a maximal
        packing set tree of depth $\overline{J}$ \emph{before} seeing the data.}
    \vspace{1mm}
    \newline
    \nit
    Construct a tree of packing sets of depth $\overline{J} \in \NN$, which is
    \emph{independent} of the data samples $\combdatavec$. Here, $\overline{J}$
    is as defined in \Cref{nthm:upper-bound-rate-finite-iterations}. The
    explicit construction of such a packing set tree proceeds as follows. First,
    fix any arbitrary point $\Upsilon_{1} \in \cF$, which is the root node, \ie,
    the first level of the packing set tree. In the case where $\overline{J} =
    1$, the tree construction stops at this single root node.

    Assuming the (more interesting) case where $\overline{J} > 1$, we then let
    $d \defined \diam{(\cF)}$, and construct a maximal $\frac{d}{2(C +
    1)}$-packing set of $B_{2}(\Upsilon_{1}, d) \cap \cF = \cF$. Denote this
    packing set by $P_{\Upsilon_{1}} \defined \theset{m_{1}, m_{2}, m_{3},
    \ldots, m_{\absb{P_{\Upsilon_{1}}}}}$. The set $P_{\Upsilon_{1}}$ forms the
    children (densities) of our root node, that is the second level of the
    tree\footnote{By \emph{convention}, the children forming the packing set
    densities are arbitrarily indexed in an increasing alphanumeric manner, from
    left child node to right child node.}. Now, for \emph{each} density in
    $P_{\Upsilon_{1}}$, we again construct a maximal packing set as follows. For
    example, taking the density $m_{3} \in P_{\Upsilon_{1}}$, we construct a
    maximal $\frac{d}{4(C + 1)}$-packing set of $B_{2}(m_{3}, d/2) \cap \cF$,
    which we denote as $P_{m_{3}} \defined \theset{m_{3,1}, m_{3,2}, m_{3,3},
    \ldots, m_{3,\absb{P_{m_{3}}}}}$. Here, the (finite) packing set $P_{m_{3}}$
    again forms the children of the node density $m_{3}$. Iterating this process
    over each density in $P_{\Upsilon_{1}}$, forms the complete second level of
    the tree.

    Now we can further iterate this process over each density in the second
    level of the tree to construct the third level of the tree. For example,
    taking the density $m_{3, 3}$, we construct a maximal $\frac{d}{8(C +
    1)}$-packing set of $B_{2}(m_{3,3}, d/4) \cap \cF$, which we denote as
    $P_{m_{3, 3}} \defined \theset{m_{3,3,1}, m_{3,3,2}, m_{3,3,3}, \ldots,
    m_{3,3\absb{P_{m_{3,3}}}}}$, which forms the children of node $m_{3, 3}$.
    This process is iterated so that for the $k^{\textnormal{th}}$-level of the
    tree, we construct $\frac{d}{2^{k}(C + 1)}$-packing sets, with closed balls
    $B_{2}(\cdot, d/2^{k - 1}) \cap \cF$. In particular the packing set tree is
    extended for each depth level $k \in \theset{2, 3, \ldots, \overline{J}-1}$.
    This process results in a maximal packing set tree of depth $\overline{J}$,
    as claimed.
    \item \textbf{Build a finite sequence of densities by traversing our packing
        set tree.}
    \label{itm:likelihood-based-estimator-03}
    \vspace{1mm}
    \newline
    \nit
    Now, \emph{after} observing our data sample $\combdatavec$, we construct a
    \emph{finite} sequence of densities, \ie, $\Upsilon \defined
    \theseqb{\Upsilon_{k}}{k = 1}{\overline{J}}$, using our packing set tree
    construction in \ref{itm:likelihood-based-estimator-02}. First, we
    initialize the first term of our sequence to $\Upsilon_{1}$, \ie, the root
    node already chosen in \ref{itm:likelihood-based-estimator-02}. If
    $\overline{J} = 1$, then the sequence $\Upsilon \defined (\Upsilon_{1})$.
    Otherwise, if $\overline{J} > 1$, we traverse down one level of our packing
    set tree, and assign $\Upsilon_{2}$ to be the density from
    $P_{\Upsilon_{1}}$ which maximizes the log-likelihood \emph{given} the data.
    That is, set $\Upsilon_{2} \defined \argmax_{\nu \in P_{\Upsilon_{1}}}
    \sum_{i = 1}^{n} \log \nu(X_{i})$. Since $P_{\Upsilon_{1}}$ is a finite set,
    this will be exhausted for each such iteration in finitely many steps.

    Moreover, we note that when assigning $\Upsilon_{2}$, there may be ties in
    children densities who all simultaneously maximize the log-likelihood. To
    break ties, by \emph{convention}, we always select the \emph{left-most}
    child from our packing set tree\footnote{This selection rule thus
    effectively assigns the child density maximizing log-likelihood with the
    \emph{smallest} such alphanumeric index.}. Once the $\Upsilon_{2}$ is
    assigned from our packing set tree, once again assign $\Upsilon_{3}$ from
    its children by again maximizing the log-likelihood. Keep iterating in this
    manner for each index\footnote{Note that $k$ here refers to index of the
    $k^{\textnormal{th}}$-term our sequence $\Upsilon$.} $k \in \theset{2, 3,
    \ldots, \overline{J}}$, and construct the finite, \ie, \emph{terminating}
    sequence $\Upsilon$.
    \item \textbf{Output estimator as the $\overline{J}^{\textnormal{th}}$-term
    of the sequence.}
    \label{itm:likelihood-based-estimator-04}
    \vspace{1mm}
    \newline
    \nit
    Finally, we note that the finite sequence $\Upsilon \defined
        \theseqb{\Upsilon_{k}}{k = 1}{\overline{J}}$ satisfies\footnote{We will
        formally justify this in the appendix in
        \Cref{nlem:upsilon-is-cauchy-sequence}.} $\normb{\Upsilon_{J} -
        \Upsilon_{J^{\prime}}}_{2}
        \leq
        \frac{d}{2^{J^{\prime} - 2}}$, for each pair of positive integers
    $J^{\prime} < J$. Our multistage sieve MLE, \ie, $\nu^{*}(\combdatavec)$,
    can be taken as the final term of this sequence. That is
    $\nu^{*}(\combdatavec) \defined \Upsilon_{\overline{J}}$. The estimator
    $\nu^{*}(\combdatavec)$ is readily understood by comparing\footnote{We note
    that in \Cref{fig:packing-set-tree-construction} if $\overline{J} = 1$, the
    estimator would just output $\Upsilon_{1}$. In the case where $\overline{J}
    > 1$, the maximal packing sets for each level of the tree are illustrated on
    the left, and the corresponding constructed tree level is shown on the
    right. In this instance the finite sequence of $\overline{J}$ densities is
    given by $\Upsilon = \parens{\Upsilon_{1}, m_{3}, m_{3, 3}, m_{3, 3, 2},
    \ldots, m_{3, 3, 2, \ldots, 5}}$. The estimator then takes the
    $\overline{J}^{\textnormal{th}}$-term of $\Upsilon$, \ie,
    $\nu^{*}(\combdatavec)=m_{\underbrace{3,3,2, \ldots,
    5}_{(\bar{J}-1)-\text{terms}}}$.} \Cref{fig:packing-set-tree-construction}
    with the qualitative description in
    \ref{itm:likelihood-based-estimator-01}-\ref{itm:likelihood-based-estimator-04}.
    \eenum
\end{mycolorbox}

\iftoggle{showfigs}{%
    \begin{figure}[!ht]
        \centering
        \includegraphics[scale=\iftoggle{cmuthesismode}{0.65}{0.70}]{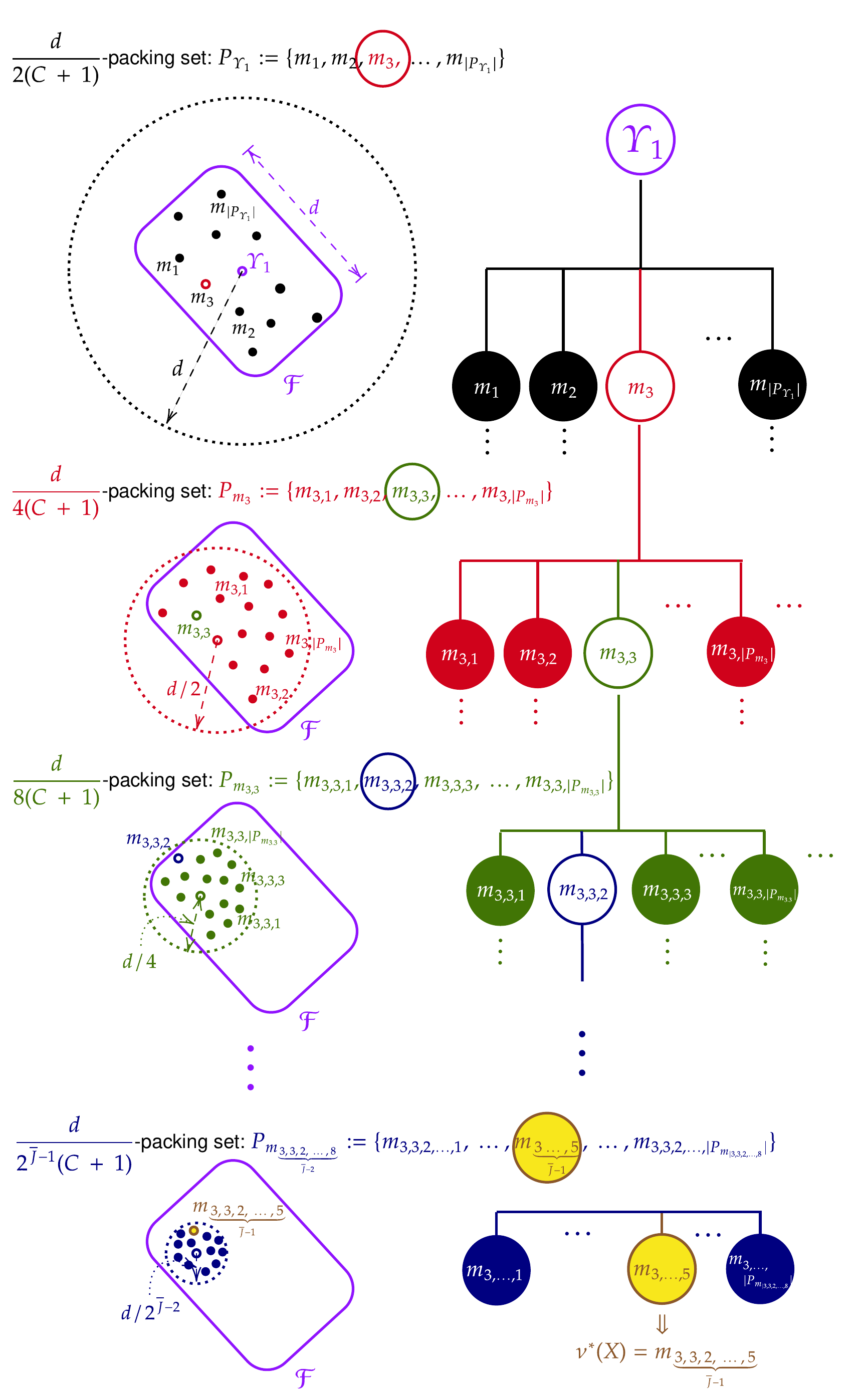}
        \caption{Maximal packing set tree construction in
            \ref{itm:likelihood-based-estimator-02}. The color of nodes
            corresponds to the color of the points. The nodes with a white
            background correspond to the points we select at each iteration of
            the multistage sieve estimation. The leaf node with the yellow
            background represents the final multistage sieve estimator.}
        \label{fig:packing-set-tree-construction}
    \end{figure}
}{%
}

\begin{restatable}{nrmk}{packingsettreeconstructionscale}\label{nrmk:packing-set-tree-construction-scale}
    We reiterate that \Cref{fig:packing-set-tree-construction} is not drawn to
    any precise scale. In reality the $L_2$-balls should be much ``wider'' than
    the set $\cF$ (and $\cFalphabetaB$). This is because they do not impose that
    their elements are proper densities, unlike the elements of the set $\cF$
    (and $\cFalphabetaB$) which are non-negative and integrate to $1$. It is
    intended to be useful conceptual guide to understanding the construction of
    our multistage sieve MLE.
\end{restatable}

We observe that our proposed estimator $\nu^{*}(\combdatavec)$ can be thought of
as an ``multistage sieve MLE'' in the spirit of \citet{wong1995probability}.
Broadly speaking a `sieve' MLE effectively takes the MLE over a strategically
constrained subset of the parameter space, \ie, $\cF$ in our setting \citep[see
Chapter~8][\eg, for more details]{grenander1981abstractinference}. Specifically,
as we traverse the down the finite-depth maximal packing set tree, each group of
children densities along with the MLE selection rule can be thought of as a
``sieve''. We note that the sieve MLE proposed in \citet{wong1995probability} is
a construction which is also not practically computable for general density
classes $\cF$.

\begin{restatable}[An \emph{online} finite packing set tree
        construction]{nrmk}{packingsettreeonline}\label{nrmk:packing-set-tree-online}
        We note that the finite-depth maximal packing set tree described in
        \ref{itm:likelihood-based-estimator-02}, can be replaced with a
        conceptually simpler \emph{online} finite-depth maximal packing set tree
        construction. This proceeds as follows. Once again, as per
        \ref{itm:likelihood-based-estimator-02}, we can initialize $\Upsilon_{1}
        \in \cF$ to be the root node independently of the data. We then
        construct the second level of our packing set tree, \ie,
        $P_{\Upsilon_{1}} \defined \theset{m_{1}, m_{2}, m_{3}, \ldots,
        m_{\absb{P_{\Upsilon_{1}}}}}$, as the previously described maximal
        packing set. This first level is constructed without looking at the data
        samples $\combdatavec$. This time however, we can traverse down the
        first level of the tree and set $\Upsilon_{2} \defined \argmax_{\nu \in
        P_{\Upsilon_{1}}} \sum_{i = 1}^{n} \log \nu(X_{i})$, \ie, by using the
        data samples $\combdatavec$. Given $\Upsilon_{2}$ selected in this data
        driven manner, we can construct the second level of the tree as the
        children of $\Upsilon_{2}$, \ie, the maximal packing set
        $P_{\Upsilon_{2}}$ \emph{without} using the data samples. We can then
        set $\Upsilon_{3} \defined \argmax_{\nu \in P_{\Upsilon_{2}}} \sum_{i =
        1}^{n} \log \nu(X_{i})$, once again using the data. We can thus repeat
        this recursive process for $\overline{J}$ iterations, whereby the
        maximal packing set of children of each parent node are constructed
        without seeing the data. The specific child node is selected after
        seeing the data, \emph{and then} the estimator can traverse to one of
        these children. This does not require the all possible children of all
        possible parent nodes of the maximal packing set tree to be constructed
        up front as described in \ref{itm:likelihood-based-estimator-02}.
        Instead, we only construct the children as required in a simple
        \emph{sequential} manner.
\end{restatable}

We next show that our multistage sieve MLE is a measurable function of the data
with respect to the Borel $\sigma$-field on $\cF$ in $L_{2}$-metric topology.
This is important, because all upper bound risk rates in expectation for
$\nu^{*}(\combdatavec)$ that follow, are with respect to the $L_{2}$-metric
topology on $\cF$.

\begin{restatable}[Measurability of
        $\nu^{*}(\combdatavec)$]{nprop}{estimatormeasurability}\label{nprop:nu-star-estimator-measurability}
        The multistage sieve MLE, \ie, $\nu^{*}(\combdatavec)$, is a measurable
        function of the data with respect to the Borel $\sigma$-field on $\cF$
        in the $L_{2}$-metric topology.
\end{restatable}

With the measurability of $\nu^{*}(\combdatavec)$ established, the main theorem
establishing the performance of $\nu^{*}(\combdatavec)$ is
\Cref{nthm:upper-bound-rate-finite-iterations} below.

\begin{restatable}[Upper bound rate for the multistage sieve MLE
        $\nu^{*}(\combdatavec)$]{nthm}{upperboundratefiniteiters}\label{nthm:upper-bound-rate-finite-iterations}
        Let, $\nu^{*}(\combdatavec) = \Upsilon_{\overline{J}}$ be the output of
        the multistage sieve MLE which is run for $\overline{J} \in \NN$ steps.
        Here $\overline{J}$ is defined as the maximal integer $J \in \NN$, such
        that $\varepsilon_J \defined \frac{\sqrt{L(\alpha, \beta, c/2 - 1)}
        d}{2^{(J-2)}c}$ satisfies\footnote{Observe that by the definition of
        $\varepsilon_{\overline J}$ and \eqref{upper:bound:suff:cond} we have
        that all packing sets used in the construction of the estimator must be
        finite, even though we are not assuming that the set $\cF$ is totally
        bounded.}
    \begin{align}\label{upper:bound:suff:cond}
        n \varepsilon_J^2
        >
        2 \log \mlocc{\cF}{\varepsilon_J
            \frac{c}{\sqrt{L(\alpha, \beta, c / 2 - 1)}}}{c}
        \vee \log 2,
    \end{align}
    or $\overline{J} = 1$ if no such $J$ exists. Then
    \begin{align*}
        \EE \normb{\nu^{*}(\combdatavec) - f}^{2}_{2}
        \leq
        \bar C \varepsilon^{*2},
    \end{align*}
    for some universal constant $\bar C$, and where $\varepsilon^* \defined
        \varepsilon_{\overline{J}}$. We remind the reader that $c \defined 2(C +
        1)$ is the constant from the definition of local metric entropy, which
        is assumed to be sufficiently large. Here $C$ is assumed to satisfy
        \eqref{neqn:log-likelihood-bernstein-02a}, and $L(\alpha, \beta, C)$ is
        defined as per \eqref{neqn:log-likelihood-bernstein-02b}.
\end{restatable}

We will now formally illustrate that the above estimator achieves the minimax
rate. The precise expression of the rate is quantified in the following result.

\begin{restatable}[Minimax
        rate]{nthm}{sharpminimaxrate}\label{nthm:sharp-minimax-rate} Define
        $\varepsilon^* \defined \sup \thesetc{\varepsilon}{n \varepsilon^{2}
        \leq \log{\mlocc{\cF}{\varepsilon}{c}}}$, where $c$ in the definition of
        local metric entropy is a  sufficiently large absolute constant. Then
        the minimax rate is given by $\varepsilon^{*2} \wedge d^2$ up to
        absolute constant factors.
\end{restatable}

\begin{restatable}[Extending results to loss functions in
        $\mathsf{KL}$-divergence and the Hellinger
        metric]{nrmk}{extendkldivhellinger}\label{nthm:extend-kldiv-hellinger}
        Recall that by \Cref{nlem:equiv-kl-euc-metric} we have the ``topological
        equivalence'' of the $\mathsf{KL}$-divergence and squared Hellinger
        metric with the squared $L_{2}$-metric on $\cF$. This means that we can
        readily extend our minimax risk bounds in \Cref{nthm:sharp-minimax-rate}
        to loss functions measured via $\mathsf{KL}$-divergence and the squared
        Hellinger metric. The important consideration is that
        \eqref{upper:bound:suff:cond} is still solved (in both cases) using the
        local metric entropy of $\cF$ using the squared $L_{2}$-metric. Note
        that for the $\mathsf{KL}$-divergence to be well-defined, we require
        that all densities are strictly positively lower bounded over the common
        compact support.
\end{restatable}

We now argue that the minimax rate for a class $\cF \subset \cFzerobetaB$ which
is convex and not necessarily lower bounded by $\alpha > 0$ is given by the same
equation, as long as there exists a single density in $f_{\alpha} \in \cF$ which
is $\alpha$-lower bounded. The argument used to establish this claim essentially
the same as used in \citet[Lemma~1]{yang1999information}, which we formalize for
our setting in \Cref{nprop:extend-zero-bounded-densities}. For completeness, we
provide all details for our setting in the Appendix.

\begin{restatable}[Extending results to
        $\cFzerobetaB$]{nprop}{extendzeroboundeddensities}\label{nprop:extend-zero-bounded-densities}
        Let $\cF \subset \cFzerobetaB$ be a convex class of densities, with at
        least one $f_{\alpha} \in \cF$ that is $\alpha$-lower bounded, with
        $\alpha > 0$. Then the minimax rate in the squared $L_{2}$-metric is
        $\varepsilon^{*2} \wedge d^{2}$, where $\varepsilon^* \defined \sup
        \thesetc{\varepsilon}{n \varepsilon^{2} \leq
        \log{\mlocc{\cF}{\varepsilon}{c}}}$.
\end{restatable}

\subsection{Adaptivity}\label{sec:adaptivity}

In this section we illustrate that the estimator, $\nu^{*}(\combdatavec)$, as
defined in \Cref{sec:upper-bound} is adaptive to the true density $f$. Before
that, similar to \citet{neykov2022minimax} we re-define the notion of
\emph{adaptive $L_{2}$-local metric entropy} for any density $\theta \in \cF$.

\begin{restatable}[Adaptive Local Entropy]{ndefn}{localpackingsets}
Let $\theta \in \cF$ be a density. Let $M(\theta, \varepsilon, c)$ denote the
maximal cardinality of a packing set of the set $B_{2}(\theta, \varepsilon) \cap
\cF$ at an $L_2$ distance $\varepsilon/c$. Let
\begin{align*}
    \madlocc{\cF}{\theta}{\varepsilon}{c}
    \defined
    M(\varepsilon/c, \cF \cap B_{2}(\theta, \varepsilon))
    \defines
    \mgloc{\cF \cap B_{2}(\theta, \varepsilon)}{\varepsilon/c}.
\end{align*}
We refer to $\log{\madlocc{\cF}{\theta}{\varepsilon}{c}}$ as the \emph{adaptive
    $L_{2}$-local metric entropy} of $\cF$ at $\theta$. We note that in contrast
    to \Cref{ndefn:local-entropy}, here we do not take the supremum over all
    $\theta \in \mclF$, \ie, $\madlocc{\cF}{\theta}{\varepsilon}{c}$ here
    depends on input $\theta$ value.
\end{restatable}

\begin{restatable}[Adaptive upper bound rate for the multistage sieve MLE
        $\nu^{*}(\combdatavec)$]{nthm}{adaptiveupperboundratefiniteiters}\label{nthm:adaptive-upper-bound-rate-finite-iterations}
        Let, $\nu^{*}(\combdatavec) = \Upsilon_{\overline{J}}$ be the output of
        the multistage sieve MLE which is run for $\overline{J}$ iterations
        where $\overline{J}$ is defined as the maximal solution to
        \begin{align*}
            n \varepsilon_J^2
        >
        2 \inf_{f \in \cF} \madlocc{\cF}{f}{2\varepsilon_J
        \frac{c}{\sqrt{L(\alpha, \beta, c / 2 - 1)}}}{2c}
        \vee \log 2,
        \end{align*}
        where $\varepsilon_J \defined \frac{\sqrt{L(\alpha, \beta, c/2 - 1)}
        d}{2^{(J-2)}c}$ and $\overline{J}= 1$ if no such $J$
        exists\footnote{Note that running the estimator with $\overline{J}$ many
        steps, may result into having non-finite packing sets --- that is not an
        issue however.}. Let $J^*$ be defined as the maximal integer $J \in
        \NN$, such that $\varepsilon_J \defined \frac{\sqrt{L(\alpha, \beta, c/2
        - 1)} d}{2^{(J-2)}c}$ such that\footnote{Observe that by the definition
        of $\varepsilon_{\overline J}$ and \eqref{upper:bound:suff:cond} we have
        that \emph{some} packing sets used in the construction of the
        \emph{adaptive} estimator may not be finite, but will be at most
        countable. This follows from the $L_{2}$-separability of $\cF$ and is
        formalized in \Cref{nlem:separability-of-class-F} in the appendix. We
        note that the measurability of the adaptive estimator still holds as per
        \Cref{nprop:nu-star-estimator-measurability} in this case.}.,
    \begin{align}\label{neqn:adaptive-upper-bound-suff-cond}
        n \varepsilon_J^2
        >
        2 \madlocc{\cF}{f}{2\varepsilon_J
        \frac{c}{\sqrt{L(\alpha, \beta, c / 2 - 1)}}}{2c}
        \vee \log 2,
    \end{align}
    and $J^* = 1$ if no such $J$ exists. Then
    \begin{align*}
        \EE \normb{\nu^{*}(\combdatavec) - f}^{2}_{2}
        \leq
        \bar C \varepsilon^{*2},
    \end{align*}
    for some universal constant $\bar C$, and where $\varepsilon^* \defined
        \varepsilon_{J^*}$. We remind the reader that $c \defined 2(C + 1)$ is
        the constant from the definition of local metric entropy, which is
        assumed to be sufficiently large. Here $C$ is assumed to satisfy
        \eqref{neqn:log-likelihood-bernstein-02a}, and $L(\alpha, \beta, C)$ is
        defined as per \eqref{neqn:log-likelihood-bernstein-02b}.
\end{restatable}

\section{Examples}\label{sec:examples}

We will now apply our work to derive risk bounds for density estimation (under
the squared $L_{2}$-metric) for various examples of convex density classes
$\cF$. To that end, per \Cref{nprop:extend-zero-bounded-densities} our risk
bounds only require us to establish that the stated class $\cF$ is indeed
convex, and importantly that there exists at least one density $f_{\alpha} \in
\cF$ that is positively bounded away from 0 over the entire support $B$. In
order to establish the latter fact we can usually take $f_{\alpha} \sim
\distUnif{[B]}$, and check that it lies in our density class $\cF$, and by
suitably expanding our ambient space\footnote{We reiterate that our use of
$\theset{\alpha, \beta}$ in this section (and throughout the paper) is
\emph{only} used to refer to the constants in
\Cref{ndefn:density-class-F-alpha-beta}, of $\cFalphabetaB$ and thus $\cF$.}
$\cFalphabetaB$. We will also use the following key fact relating
$L_{2}$-\emph{local} and $L_{2}$-\emph{global} metric entropies.
\begin{equation}\label{neqn:local-global-entropy-bounds}
    \log{\mgloc{\cF}{\varepsilon / c}} - \log{\mgloc{\cF}{\varepsilon}}
    \leq
    \log{\mlocc{\cF}{\varepsilon}{c}}
    \leq
    \log{\mgloc{\cF}{\varepsilon / c}}
\end{equation}
Here, \eqref{neqn:local-global-entropy-bounds} follows directly from
\citet[Lemma~2]{yang1999information}, where it is only proved for the case $c =
2$. However, their proof directly extends to the more general case for each $c >
0$, which is required for our setting. For the various examples of $\cF$ that
follow below, we will formally show\footnote{See \Cref{appendix:proofs-examples}
for these details.} that a stronger sufficient condition on global entropy is
satisfied, namely
\begin{equation}\label{neqn:local-global-entropy-bounds-suff-condn}
    \log{\mgloc{\cF}{\varepsilon / c}} - \log{\mgloc{\cF}{\varepsilon}}
    \asymp
    \log{\mgloc{\cF}{\varepsilon / c}},
\end{equation}
provided we take $c$ to be sufficiently large enough, which is within our
control to do, per our packing set tree construction. In short,
\eqref{neqn:local-global-entropy-bounds-suff-condn} will enable us to bound the
local metric entropy via \eqref{neqn:local-global-entropy-bounds}. To illustrate
this, we initially consider two examples from
\citet[see][Section~6]{yang1999information}. We begin with the class $\cF
\defined \lipschitzpsi$, \ie, the $(\gamma, q, \Psi)$-Lipschitz density class
defined as per \eqref{neqn:lipschitz-density-class-01}. As noted in
\citet[Section~6.4]{yang1999information}, with fixed constants $\max{\theset{1 /
q - 1 / 2, 0}} < \gamma \leq 1$, and $1 \leq q \leq \infty$, the
$\varepsilon$-\emph{global} metric entropy of $\lipschitzpsi$ is of the order
$\varepsilon^{-1 / \gamma}$ per \citet{birman1980quantsobolevimbedding}.
\begin{restatable}[Lipschitz density class
        $\cF$]{nexa}{exalipschitzdensityclass}\label{nexa:lipschitz-density-class}
        Let $1 < \Psi < \beta < \infty$, $\max{\theset{1 / q - 1 / 2, 0}} <
        \gamma \leq 1$, and $1 \leq q \leq \infty$ be fixed constants, and $B
        \defined [0, 1]$. Now, let $\cF \defined \lipschitzpsi$ denote the space
        of $(\gamma, q, \Psi)$-Lipschitz densities with total variation at most
        $\beta$. That is,
    \begin{equation}\label{neqn:lipschitz-density-class-01}
        \lipschitzpsi
        \defined
        \thesetb{f \colon B \to [0, \Psi]}{%
            \normb{f(x + h) - f(x)}_{q} \leq \Psi h^{\gamma},
            \normb{f}_{q} \leq \Psi,
            \int_{B} f \dlett{\mu} = 1,
            f \text{ measurable}},
    \end{equation}
    and $\normb{f}_{q} \defined \parens{\int_{B} \absa{f(x)}^{q} \dlett{\mu}}^{1
            / q}$. Note that in \eqref{neqn:lipschitz-density-class-01} we have
            that $x \in B$, and only consider $h > 0$, and further $f(x + h) =
            f(1)$, for $x + h > 1$, so that the predicate of $\lipschitzpsi$ is
            well-defined. Then $\lipschitzpsi$ is a convex density class, there
            exists a density $f_{\alpha} \in \lipschitzpsi$ that is strictly
            positively bounded away from 0, and the minimax rate (in the squared
            $L_{2}$-metric) for estimating $f \in \lipschitzpsi$ is of the order
            $n^{- \frac{2 \gamma}{2 \gamma + 1}}$.
\end{restatable}

Another well studied density estimation problem is the case where $\cF \defined
    \totbddvzeta$ is total bounded variation at most $\zeta$, defined as per
    \eqref{neqn:bdd-total-variation-density-class-01}. Importantly we note that
    the $\varepsilon$-\emph{global} $L_{2}$-metric entropy of this well studied
    function class is of the order $\varepsilon^{-1}$ \citep[see
    Section~6.4][\eg, for more details]{yang1999information}.
\begin{restatable}[Bounded total variation density class
        $\cF$]{nexa}{exabddtotvardensityclass}\label{nexa:bdd-total-variation-density-class}
        Let $1 < \zeta < \beta < \infty$ be a fixed constant, and $B \defined
        [0, 1]$. Now, let $\cF \defined \totbddvzeta$ denote the space of
        univariate densities with total variation at most $\beta$. That is,
    \begin{equation}\label{neqn:bdd-total-variation-density-class-01}
        \totbddvzeta
        \defined
        \thesetb{f \colon B \to [0, \zeta]}{%
            \normb{f}_{\infty} \leq \zeta,
            V(f) \leq \zeta,
            \int_{B} f \dlett{\mu} = 1,
            f \text{ measurable}},
    \end{equation}
    where we define the total variation of $f$, \ie, $V(f)$ as
    \begin{equation}\label{neqn:bdd-total-variation-density-class-02}
        V(f)
        \defined
        \sup_{\thesetb{x_{1}, \ldots, x_{m}}{0 \leq x_1 < \cdots < x_m \leq 1, m \in \NN}}
        \sum_{i=1}^{m - 1} \absb{f\parens{x_{i+1}}-f\parens{x_{i}}},
    \end{equation}
    and $\normb{f}_{\infty} \defined \sup_{x \in B} \absa{f(x)}$. Then the
    minimax rate  (in the squared $L_{2}$-metric) for estimating $f \in
    \totbddvzeta$ is of the  order $n^{- 2 / 3}$.
\end{restatable}

Another interesting example illustrating the use case of our bounds is that
where $\cF \defined \quadfuncclass$, forms the density class of
$\gamma$-quadratic functionals defined as per
\eqref{neqn:quad-functional-density-class-01}. Importantly we note that the
$\varepsilon$-\emph{global} $L_{2}$-metric entropy of this well studied function
class is of the order $\varepsilon^{-1 / 4}$ \citep[see Example~15.8 and
Example~15.22][\eg, for more details]{wainwright2019high}.
\begin{restatable}[Quadratic functional density class
        $\cF$]{nexa}{quadfunctionaldensityclass}\label{nexa:quad-functional-density-class}
        Let $0 < \alpha < 1 < \beta < \infty$, and $\gamma > 1$ be fixed
        constants, with $B \defined [0, 1]$. Now, let $\cF \defined
        \quadfuncclass$ denote the space of univariate quadratic functional
        densities. That is,
    \begin{equation}\label{neqn:quad-functional-density-class-01}
        \quadfuncclass
        \defined
        \thesetb{f \colon B \to [\alpha, \beta]}{%
            \normb{f^{\prime \prime}}_{\infty} \leq \gamma,
            \int_{B} f \dlett{\mu} = 1,
            f \text{ measurable}}.
    \end{equation}
    Then $\quadfuncclass$ is a convex density class, there exists a density
    $f_{\alpha} \in \quadfuncclass$ that is strictly positively bounded away
    from 0, and the minimax rate (in the squared $L_{2}$-metric) for estimating
    $f \in \quadfuncclass$ is of the order $n^{- 4 / 5}$.
\end{restatable}

We now turn our attention to an interesting example, which demonstrates that our
results can yield useful bounds in cases where $L_{2}$-\emph{global} metric
entropy of $\cF$ may be unknown (or difficult to compute), but the
$L_{2}$-\emph{local} metric entropy can be controlled.
\begin{restatable}[Convex mixture density class
        $\cF$]{nexa}{exaconvexmixturedensityclass}\label{nexa:convex-mixture-density-class}
        Let $\cF \defined \convclassk$ where
    \begin{equation}\label{neqn:convex-mixture-density-class-01}
        \convclassk
        \defined
        \thesetb{\sum_{i = 1}^k \alpha_{i} f_{i}}{\sum_{i = 1}^k \alpha_{i} = 1,
            \alpha_{i} \geq 0,
            f_{i} \in  \cFalphabetaB},
    \end{equation}
    for some \underline{fixed} $k \in \NN$ and $f_i \in \cFalphabetaB$ for each
    $i \in [k]$. Further, let $\bfG = \parens{\bfG_{ij}}_{i, j \in [k]}$ denote
    the Gram matrix with $\bfG_{ij} \defined \int_B f_{i} f_{j} \mu(\dlett{x})$,
    which we assume is positive definite, \ie, $\bfG \succ \bfzero$. Then the
    minimax rate for estimating $f \in \convclassk$ is bounded from above by
    $\sqrt{\frac{k}{n}}$ up to absolute constant factors.
\end{restatable}

\section{Discussion}\label{sec:discussion}

In this paper we derived exact minimax rates for density estimation over convex
density classes. Our work builds on seminal research of
\citet{lecam1973convergence,birge1983approximation,yang1999information,wong1995probability}.
More directly, we non-trivially adapted the techniques of
\cite{neykov2022minimax}, who used it for deriving exact rates for the Gaussian
sequence model. Our results demonstrate that the $L_{2}$-\emph{local} metric
entropy \emph{always} determines that minimax rate under squared $L_{2}$-loss in
this setting. We thus provide a unifying perspective across parametric
\emph{and} nonparametric convex density classes, and under weaker assumptions
than those used by \citet{yang1999information}.

An important open question that we would like to think further about is whether
there exists a computationally tractable estimator which is also minimax optimal
in our setting. We can also consider applying our techniques to the
nonparametric regression setting (with Gaussian noise) where $f$ is a uniformly
bounded regression function of interest. We leave these exciting directions for
future work. Finally, we hope that this research stimulates further activity in
approximating $L_{2}$-\emph{local} metric entropy for various convex density
classes.


\section{Acknowledgments}\label{sec:acknowledgments-density-estimation-convex-class}

We would like to thank Arun Kumar Kuchibhotla for suggesting to us to prove that
our estimator is adaptive to the true density. We would also like to thank
Wanshan Li, Yang Ning, Alex Reinhart, Alessandro Rinaldo, and Larry Wasserman
for providing insightful feedback and encouragement during this work. All
figures in this paper were drawn using the
\texttt{Mathcha}\footnote{\url{https://www.mathcha.io/editor}} editor.

\clearpage
\appendix
\section{Preliminary}\label{app:mathematical-preliminaries}

We begin with some basic mathematical preliminaries for our work.

\iftoggle{cmuthesismode}{%
 \subsection{Notation
 Summary}\label{subsec:app-notation-summary-density-estimation-convex-class}

 To ensure that the Appendix is can be read in a standalone manner, we
 consolidate key notation used in the paper in
 \Cref{tab:notation-density-estimation-convex-class}.

 \begin{table}[htbp]\caption{Notation and conventions used in
 \iftoggle{cmuthesismode}{this chapter}{our paper}} \centering
     \begin{tabular}{r l p{10cm} }
         \toprule
         \multicolumn{2}{c}{\underline{\textbf{Variables and inequalities}}} \\
         \multicolumn{2}{c}{}
         \\
         $a \wedge b$                                       & $\min\braces{a,
         b}$ for each $a, b \in \reals$ \\
         $a \vee b$                                         & $\max\braces{a,
         b}$ for each $a, b \in \reals$ \\
         $\lesssim$                                         & $\leq$ up to
         positive universal constants \\
         $\gtrsim$                                          & $\geq$ up to
         positive universal constants \\
         $\asymp$                                          & if both $\lesssim$
         and $\gtrsim$ hold \\
         \multicolumn{2}{c}{}
         \\
         \multicolumn{2}{c}{\underline{\textbf{Functions and sets}}} \\
         \multicolumn{2}{c}{}
         \\
         $\normb{\cdot}_{2}$                                & the $L_{2}$-metric
         in $\cF$ \\
         $[m]$                                              & $\theseta{1,
         \ldots, m}$, for $m \in \nats$ \\
         $B_{2}(\theta, r)$                                 & closed
         $L_{2}$-ball centered at $\theta \in \cF$ with radius $r$ \\
         \bottomrule
     \end{tabular}
     \label{tab:notation-density-estimation-convex-class}
 \end{table}
}{%
}

\subsection{Properties of
    \texorpdfstring{$\cFalphabetaB$}{cFalphabetaB}}\label{app:prop-of-class-F-alpha-beta}
    Here we provide some basic analytic properties of our core density class
    $\cFalphabetaB$, as per \Cref{ndefn:density-class-F-alpha-beta}. Many of
    these facts will be used, sometimes implicitly, in our proofs. We hope that
    by documenting them rigorously, they provide the reader with a much richer
    understanding of the geometry of this broader density class. This may also
    be a useful reference for researchers in working in similar density
    estimation settings. We provide suitable references where the properties
    follow from standard real analysis theory.

\begin{restatable}[Convexity of
        $\cFalphabetaB$]{nlem}{convexitycfalphabeta}\label{nlem:convexity-of-class-F-alpha-beta}
        The density class $\cFalphabetaB$, forms a convex set, in the
        $L_{2}$-metric.
\end{restatable}
\bprfof{\Cref{nlem:convexity-of-class-F-alpha-beta}} In order to show the
convexity of $\cFalphabetaB$, Let $f, g \in \cFalphabetaB$, and let $\kappa \in
[0, 1]$ be arbitrary. Then for each $x \in B$, we observe that
\begin{align}
    (\kappa f + (1 - \kappa) g)(x)
    \defined \kappa f(x) + (1 - \kappa) g(x)
     & \geq \kappa \alpha + (1 - \kappa) \alpha
    \geq \alpha
    \label{neqn:prop-of-class-F-alpha-beta-01}  \\
    (\kappa f + (1 - \kappa) g)(x)
    \defined \kappa f(x) + (1 - \kappa) g(x)
     & \leq \kappa \beta + (1 - \kappa) \beta
    \leq \beta
    \label{neqn:prop-of-class-F-alpha-beta-02}
\end{align}
From \eqref{neqn:prop-of-class-F-alpha-beta-01} and
\eqref{neqn:prop-of-class-F-alpha-beta-02}, it follows that $\kappa f + (1 -
\kappa) g \colon B \to [\alpha, \beta]$. Moreover, since $\int_{B} f \dlett{\mu}
= \int_{B} g \dlett{\mu} = 1$, we have
\begin{equation}
    \int_{B} (\kappa f + (1 - \kappa) g) \dlett{\mu}
    = \kappa \int_{B} f \dlett{\mu} +
    (1 - \kappa) \int_{B} g \dlett{\mu}
    = 1.
\end{equation}
Since $f, g$ are measurable functions, then so is their convex combination, \ie,
$\kappa f + (1 - \kappa) g$. Combining the above we have shown that $\kappa f +
(1 - \kappa) g \in \cFalphabetaB$, which proves the convexity of
$\cFalphabetaB$, as required.
\eprfof

\begin{restatable}[Boundedness of
        $\cFalphabetaB$]{nlem}{boundednesscfalphabeta}\label{nlem:boundedness-of-class-F-alpha-beta}
        The density class $\cFalphabetaB$, is bounded, in the $L_{2}$-metric.
\end{restatable}
\bprfof{\Cref{nlem:boundedness-of-class-F-alpha-beta}} We now show that
$\cFalphabetaB$ is bounded in the $L_{2}$-metric. To see this observe that for
any $f, g \in \cFalphabetaB$:
\begin{equation}
    \normb{f - g}_{2}^{2}
    \defined \int_{B} (f-g)^{2} \dlett{\mu}
    \leq \int_{B} \absa{f-g} 2 \beta \dlett{\mu}
    \leq 2 \beta \parens{\int_{B} \absa{f} \dlett{\mu} + \int_{B} \absa{g} \dlett{\mu}}
    = 4 \beta.
\end{equation}
It follows that $\diam_{2}\parens{\cFalphabetaB} \defined \sup\thesetb{\normb{f
            - g}_{2}}{f, g \in \cFalphabetaB} \leq 2 \sqrt{\beta} < \infty$, as
            required.
\eprfof

\begin{restatable}[$\cFalphabetaB$ lies in
        $L^{2}(B)$]{nlem}{subsetlebtwospacecfalphabeta}\label{nlem:subset-l2-space-class-F-alpha-beta}
        The density class $\cFalphabetaB$, satisfies $\cFalphabetaB \subset
        L^{2}(B)$, where
    \begin{equation}\label{neqn:l2-space-on-set-b-01}
        L^{2}(B)
        \defined
        \thesetb{f \colon B \to \RR}{\int_{B}{f}^{2} \dlett{\mu} < \infty, \textnormal{$f$ measurable}}.
    \end{equation}
    As such $(\cFalphabetaB, \normb{\cdot}_{2})$ is an induced metric subspace
    of $(L^{2}(B), \normb{\cdot}_{2})$.
\end{restatable}
\bprfof{\Cref{nlem:subset-l2-space-class-F-alpha-beta}} Let $f \in
    \cFalphabetaB$ be arbitrary. We then observe that
\begin{equation}
    \int_{B}{f}^{2} \dlett{\mu}
    \leq \int_{B}f {\beta} \dlett{\mu}
    = \beta
    < \infty,
\end{equation}
since $f \leq \beta$ by definition of $\cFalphabetaB$. Given that $f \colon B
    \to [\alpha, \beta] \subset \RR$, we have that $f \in L^{2}(B)$, \ie,
    $\cFalphabetaB \subset L^{2}(B)$ as required.
\eprfof

\begin{restatable}[Completeness and separability of
        $L^{2}(B)$]{nlem}{completenesslebtwospace}\label{nlem:completeness-separability-l2-space}
        The metric space $(L^{2}(B), \normb{\cdot}_{2})$, with $L^{2}(B)$
        defined as per \Cref{neqn:l2-space-on-set-b-01}, is complete and
        separable.
\end{restatable}
\bprfof{\Cref{nlem:completeness-separability-l2-space}} We note that
completeness of $(L^{2}(B), \normb{\cdot}_{2})$ follows directly from
\citet[Theorem~4.8]{brezis2011funcanalysis}, and separability follows from
\citet[Theorem~4.13]{brezis2011funcanalysis}.
\eprfof

\begin{restatable}[Completeness and separability of $(\cFalphabetaB,
            \normb{\cdot}_{2})$]{nlem}{completenesscfalphabeta}\label{nlem:completeness-of-class-F-alpha-beta}
            The metric space $(\cFalphabetaB, \normb{\cdot}_{2})$ is complete
            and separable.
\end{restatable}
\bprfof{\Cref{nlem:boundedness-of-class-F-alpha-beta}} Firstly we note that
$(\cFalphabetaB, \normb{\cdot}_{2})$ is an induced metric subspace of
$(L^{2}(B), \normb{\cdot}_{2})$ per
\Cref{nlem:subset-l2-space-class-F-alpha-beta}. Now separability of
$(\cFalphabetaB, \normb{\cdot}_{2})$ follows, since it is inherited from
$(L^{2}(B), \normb{\cdot}_{2})$ by applying
\citet[Proposition~2.3.16]{shirali2006metricspaces}. We now show the
completeness of $(\cFalphabetaB, \normb{\cdot}_{2})$. Take an arbitrary Cauchy
sequence $\theseqb{f_k}{k = 1}{\infty}$ in $\cFalphabetaB$. Since $(L^{2}(B),
\normb{\cdot}_{2})$ is complete per
\Cref{nlem:completeness-separability-l2-space}, it follows that the $L_{2}$
limit of $\theseqb{f_k}{k = 1}{\infty}$ exists in $L^{2}(B)$. Let $f$ be that
limit, \ie, $\lim_{k \to \infty} f_{k} \defines f \in L^{2}(B)$. We will show
that $f \in \cFalphabetaB$. First let us show that it is a density, \ie, it
integrates to 1. By Cauchy-Schwartz $\int_{B} |f_k(x) - f(x)| \mu(\dlett{x})
\leq \sqrt{\int_{B} \absb{f_k(x) - f(x)}^2 \mu(\dlett{x})} \rightarrow 0$ so
that $\int f(x) \mu(\dlett{x}) = 1$. Next consider the function $f^{\prime} = (f
\wedge \alpha) \vee \beta$. Since for any $x \in [\alpha ,\beta]$ and any $y$ we
have $|x - y| \geq |x - (y \wedge \alpha) \vee \beta|$ then that implies (since
$f_k \in \cFalphabetaB$) that $\int_{B} |f_k(x) - f^{\prime}(x)|^2
\mu(\dlett{x}) \leq \int_{B} |f_k(x) - f(x)|^2 \mu(\dlett{x}) \rightarrow 0$ and
so $f^{\prime}$ must also be a limit of $f_{k}$. Since the limits are unique (up
to considering equivalence classes modulo sets of measure $0$ with respect to
$\mu$) then $f = f^{\prime}$ and hence it belongs to $\cFalphabetaB$.
\eprfof

\begin{restatable}[Separability of  $(\cF,
    \normb{\cdot}_{2})$]{nlem}{completenesscf}\label{nlem:separability-of-class-F}
    The metric space $(\cF, \normb{\cdot}_{2})$ is separable. Furthermore, if $A
    \subset \cF$, then $(A, \normb{\cdot}_{2})$ is separable.
\end{restatable}
\bprfof{\Cref{nlem:separability-of-class-F}} We first observe by
\Cref{nlem:completeness-of-class-F-alpha-beta} that the metric space
$(\cFalphabetaB, \normb{\cdot}_{2})$ is separable. It then follows that since
$\cF \subset \cFalphabetaB$, that $(\cF, \normb{\cdot}_{2})$ is a restriction of
$(\cFalphabetaB, \normb{\cdot}_{2})$, and thus a separable metric space by
\citet[Proposition~2.3.16]{shirali2006metricspaces}. By a similar argument, it
follows that if $A \subset \cF$, then $(A, \normb{\cdot}_{2})$ is separable.
\eprfof

The following lemma shows that in general $\cFalphabetaB$ is not totally
bounded. We consider a restricted case of $B \defined [0, 1]$, to construct a
suitable counterexample.
\begin{restatable}[Non total boundedness of $(\cFalphabetasupp{[0, 1]},
            \normb{\cdot}_{2})$]{nlem}{nontotbddlebtwospace}\label{nlem:non-total-bdd-l2-zero-one}
            Suppose $\beta \geq 2 - \alpha$. Then the metric space
            $(\cFalphabetasupp{[0, 1]}, \normb{\cdot}_{2})$ is not totally
            bounded, and hence not compact.
\end{restatable}
\bprfof{\Cref{nlem:non-total-bdd-l2-zero-one}} We note that per
\citet[Theorem~5.1.12]{shirali2006metricspaces} a metric space is totally
bounded if and only if every sequence contains a Cauchy subsequence. We will use
this characterization to construct a counterexample to demonstrate that
$(\cFalphabetasupp{[0, 1]}, \normb{\cdot}_{2})$ is not totally bounded. In
particular, we will define a sequence in $(\cFalphabetasupp{[0, 1]},
\normb{\cdot}_{2})$ which can't contain \emph{any} Cauchy subsequence.

Specifically we consider the sequence of functions $(1, \{x\mapsto \sin(2\pi j
    x)\}_{j \in \NN})$. These functions are orthonormal in $L^2([0, 1])$.
    Construct the sequence of functions $f_{j}(x) = 1 + (1-\alpha)\sin(2\pi j
    x)$ for $j \in \NN$. By the orthogonality of $1$ and $\sin(2\pi j x)$ we
    have that $\int_0^1 f_{j}(x) dx = 1$. Furthermore, $\alpha \leq f_{j}(x)
    \leq 2 - \alpha$ for all $x \in [0,1]$, hence since $\beta \geq 2 -\alpha$
    we have $f_{j}(x) \in \mathcal{F}^{[\alpha,\beta]}_{[0,1]}$. Take any two $j
    \neq k \in \NN$, and consider
\begin{align*}
    \|f_{j} - f_k\|_2^2 = (1-\alpha)^2\|\sin(2\pi j x) - \sin(2\pi k x)\|_2^2 = 2(1-\alpha)^2 > 0.
\end{align*}
This shows that there cannot be a Cauchy subsequence and hence the set is not
totally bounded.
\eprfof

\subsection{Elementary inequalities}\label{app:elementary-inequalities}

We will state and prove \Cref{nlem:elem-log-inequality}, which will provide the
key fact to will assist us in the proof of the lower bound in
\Cref{nlem:equiv-kl-euc-metric}.

\begin{restatable}[Elementary $\log$
        inequality]{nlem}{elemlogineq}\label{nlem:elem-log-inequality} For each
        $\gamma > 0$, and for any $x \in (0, \gamma]$, the following
        relationship holds:
    \begin{equation}\label{neqn:nlem:elem-log-inequality-01}
        \log{x}
        \leq
        (x - 1) - h(\gamma) (x - 1)^{2}.
    \end{equation}
    Here $h : (0, \infty) \to \RR$ is defined as in
    \eqref{neqn:nlem:equiv-kl-euc-metric-02}, and is positive over its entire
    support.
\end{restatable}
\bprfof{\Cref{nlem:elem-log-inequality}} We first argue that $h(x) > 0$ for $x
    \in (0,\infty)$. This is by the elementary inequality $\log(x + 1) \leq x$
    for all $x \geq -1$. Next, it suffices to show that the map $x \mapsto h(x)$
    is decreasing for $x > 0$ where $h$ is defined in
    \eqref{neqn:nlem:equiv-kl-euc-metric-02}. This is because
    \eqref{neqn:nlem:elem-log-inequality-01} holds for $x = 1$, and if $ x \neq
    1$ it is equivalent to
\begin{align*}
    h(\gamma)
    \leq
    \frac{(x-1) - \log x}{(x-1)^2},
\end{align*}
for $x \leq \gamma$. It is simple to verify that
\begin{align*}
    h'(x) = \frac{-x^2 + 2x \log x + 1}{(x-1)^3 x}.
\end{align*}
We will show that the above function is negative on $(0,\infty)$ which will
complete the proof. First we will evaluate it at $x = 1$. By a triple
application of L'H\^opital's rule it is simple to verify that $\frac{d}{dx}h(x)
\vert_{x = 1} = -\frac{1}{3} < 0$. Thus, it remains to show that for $x\neq 1$,
\begin{align*}
    (-x^2 + 2x \log x + 1)(x-1) < 0.
\end{align*}
Now let $f(x) \defined \frac{x^{2}-1}{2 x} - \log{x}$. We want to show that
$f(x) > 0$, for each $x > 1$ and $f(x) < 0$ for $x < 1$. First observe that
$f(1) = 0$. Moreover, we have that
\begin{equation}
    f^{\prime}(x) = \frac{(x - 1)^{2}}{2 x^{2}} > 0.
\end{equation}
That is, $f(x)$ is \emph{strictly} increasing, which implies that $f(x) > 0$ for
each $x \in (1, \infty)$, and $f(x) < 0$ for each $x < 1$ as required.
\eprfof

\clearpage
\newpage
\section{Proofs of
\texorpdfstring{\Cref{sec:minimax-upper-and-lower-bounds}}{\autoref{sec:minimax-upper-and-lower-bounds}}}\label{appendix:proofs-minimax-upper-and-lower-bounds}

\subsection{Proof of
    \texorpdfstring{\Cref{nlem:equiv-kl-euc-metric}}{\autoref{nlem:equiv-kl-euc-metric}}}\label{app:nlem:equiv-kl-euc-metric}
\equivkleucmetric*

\bprfof{\Cref{nlem:equiv-kl-euc-metric}} We will prove the upper and lower bound
in turn.\\
\newline
\textbf{(Upper bound in \eqref{neqn:nlem:equiv-kl-euc-metric-01}):} We seek to
show that $\kldiva{f}{g} \leq \frac{1}{\alpha} \normb{f - g}_{2}^{2}$. First,
for any two densities $f, g \in \cF$, we define the $\chi^{2}$-divergence
between $f$ and $g$, as follows:
\begin{equation}\label{neqn:chi-sq-divergence-density}
    \chisqdiva{f}{g}
    \defined \int_{B} \frac{(f - g)^{2}}{g} \dlett{\mu}.
\end{equation}
Per \Cref{nrmk:kldiv-density-well-defined}, we note that $\chi^{2}(f || g)$ in
\eqref{neqn:chi-sq-divergence-density} is similarly well-defined. We then have:
\begin{align}
    \kldiva{f}{g}
     & \leq \chisqdiva{f}{g}
    \tag{per \citet[Theorem~5]{gibbs2002chooseboundprobmetrics}} \\
     & \defined \int_{B} \frac{(f - g)^{2}}{g} \dlett{\mu}
    \tag{using \eqref{neqn:chi-sq-divergence-density}}           \\
     & \leq \frac{1}{\alpha} \int_{B} (f - g)^{2} \dlett{\mu}
    \tag{since $\inf_{x \in B} g(x) \geq \alpha > 0$}            \\
     & \defines \frac{1}{\alpha} \normb{f - g}_{2}^{2}.
    \tag{by definition}
\end{align}
As required. \qedbsquare \\
\newline
\noindent\textbf{(Lower bound in \eqref{neqn:nlem:equiv-kl-euc-metric-01}):} We
seek to show that $\kldiva{f}{g} \geq c(\alpha,\beta) \normb{f - g}_{2}^{2}$. We
proceed as follows: First observe that for any $f, g \in \cF$, we have that $0 <
\frac{g}{f} \leq \frac{\beta}{\alpha} < \infty$
\begin{align}
    \kldiva{f}{g}
     & \defined
    \int_{B} f \log\parens{\frac{f}{g}} \dlett{\mu}
    \tag{per \eqref{neqn:kldiv-density-defn}}                                                          \\
     & =
    \int_{B} - f \log\parens{\frac{g}{f}} \dlett{\mu}
    \tag{since $\inf_{x \in B} f(x) \geq \alpha > 0$}                                                  \\
     & \geq
    \int_{B} - f \parens{\parens{\frac{g}{f} - 1} - h(\beta / \alpha)\parens{\frac{g}{f} - 1}^{2}} \dlett{\mu}
    \tag{using \Cref{nlem:elem-log-inequality}, with $C = \frac{\beta}{\alpha}$ and $x = \frac{g}{f}$} \\
     & =
    \int_{B} (f - g) \dlett{\mu} +
    h(\beta / \alpha) \int_{B} \frac{(g - f)^{2}}{f} \dlett{\mu}
    \nonumber                                                                                          \\
     & \geq \frac{h(\beta / \alpha)}{\beta} \int_{B} (g - f)^{2} \dlett{\mu}
    \tag{since $\int_{B} (f - g) \dlett{\mu} = 0$, and $0 < \sup_{x \in B} f(x) \leq \beta$}           \\
     & \defines \frac{h(\beta / \alpha)}{\beta} \normb{f - g}_{2}^{2}
    \nonumber                                                                                          \\
     & \defines c(\alpha, \beta) \normb{f - g}_{2}^{2},
\end{align}
where we define $c(\alpha, \beta) \defined \frac{h(\beta / \alpha)}{\beta} > 0$.
This proves the lower bound in \eqref{neqn:nlem:equiv-kl-euc-metric-01}, as
required.  \qedbsquare

We now show the following equivalence between the Hellinger, \ie,
$d_{\mathsf{H}}$-metric, and the $L_{2}$ metric in $\cFalphabetaB$.
\begin{equation}\label{neqn:equiv-hellinger-l2-metric-01}
    (1 / 4 \beta) \normb{f-g}_{2}^{2}
    \leq
    \hellmet{f}{g}^{2}
    \leq
    (1 / \alpha) \normb{f - g}_{2}^{2}.
\end{equation}
To prove the upper bound in \eqref{neqn:equiv-hellinger-l2-metric-01}, we note
that
\begin{align}
    \hellmet{f}{g}^{2}
     & \leq
    \kldiva{f}{g}
    \tag{from \citet{gibbs2002chooseboundprobmetrics}} \\
     & \leq
    (1 / \alpha) \normb{f - g}_{2}^{2},
    \tag{per \eqref{neqn:nlem:equiv-kl-euc-metric-01}}
\end{align}
as required. In order to prove the lower bound we observe that
\begin{align}
    \normb{f - g}_{2}^{2}
     & =
    \int_{B} (f - g)^{2} \dlett{\mu}
    \nonumber                      \\
     & =
    \int_{B} (\sqrt{f} + \sqrt{g})^{2} (\sqrt{f} - \sqrt{g})^{2} \dlett{\mu}
    \nonumber                      \\
     & \leq
    4 \beta \int_{B} (\sqrt{f} - \sqrt{g})^{2} \dlett{\mu}
    \tag{since $f, g \leq \beta$.} \\
     & \defines
    4 \beta \hellmet{f}{g}^{2},
    \tag{by definition}
\end{align}
which implies the required lower bound in
\eqref{neqn:equiv-hellinger-l2-metric-01}. We have thus established the required
upper and lower bounds in both \eqref{neqn:nlem:equiv-kl-euc-metric-01} and
\eqref{neqn:equiv-hellinger-l2-metric-01}.

Finally, we note that $(\cF, \normb{\cdot}_{2})$ is metric space, since it is
the restriction of the metric space $(\cFalphabetaB, \normb{\cdot}_{2})$. And so
the bounds \eqref{neqn:nlem:equiv-kl-euc-metric-01} and
\eqref{neqn:equiv-hellinger-l2-metric-01}, are also inherited by $\cF \subset
\cFalphabetaB$.
\eprfof

\subsection{Proof of
    \texorpdfstring{\Cref{nlem:minimax-lower-bound}}{\autoref{nlem:minimax-lower-bound}}}\label{app:nlem:minimax-lower-bound}

\minimaxlowerbound*
\bprfof{\Cref{nlem:minimax-lower-bound}} Let $c > 0$ be fixed, and $\theta \in
    \cF$ be an arbitrary point. Consider maximal packing the set $\theset{f^{1},
    \ldots, f^{m}} \subset \cF \cap B_{2}(\theta, \varepsilon)$ at a
    $L_{2}$-``distance'' at least $\varepsilon/c$. Here $B_{2}(\theta,
    \varepsilon)$ denotes a closed $L_{2}$-ball around the point $\theta$, with
    radius $\varepsilon$. Suppose it has $m$ elements. Then we know that
\begin{equation*}
    I(X;J)
    \leq
    \frac{1}{m}\sum_{j = 1}^{m} \kldiva{f^{j}}{\theta}
    \leq
    \max_{j \in [m]} \kldiva{f^{j}}{\theta}
    \leq
    \max_{j \in [m]} (1 / \alpha) \normb{f^{j} - \theta}_{2}^{2}
    \leq
    \varepsilon^{2} / \alpha.
\end{equation*}
Here the final two inequalities follow by applying
\eqref{neqn:nlem:equiv-kl-euc-metric-01}, and using the fact that
$\theset{f^{1}, \ldots, f^{m}} \subset \cF \cap B_{2}(\theta, \varepsilon)$,
respectively. Hence, if the packing number satisfies $\log m \geq 2 n
\varepsilon^{2} / \alpha + 2\log 2$ we will have a lower bound proportional to
$\varepsilon^2$ (it will be $\varepsilon^2/(8c^2)$). By taking the supremum over
$\theta$, we conclude that if $\log{\mlocc{\cF}{\varepsilon}{c}} > 2 n
\varepsilon^{2} / \alpha + 2\log 2$ we have a lower bound proportional to
$\varepsilon^2$.
\eprfof

\subsection{Proof of
    \texorpdfstring{\Cref{nlem:log-likelihood-bernstein}}{\autoref{nlem:log-likelihood-bernstein}}}\label{app:nlem:log-likelihood-bernstein}

\loglikelihoodhoeffding*
\bprfof{\Cref{nlem:log-likelihood-bernstein}} We first observe per
\Cref{nrmk:log-likelihood-well-defined} that the log-likelihood, $\psi(g,
g^{\prime}, X)$, is well-defined. Next the mean of these variables, for each $i
\in [n]$, is
\begin{align}
    \EE_f \brackets{\log \frac{g(X_{i})}{g^{\prime}(X_{i})}}
     & = \EE_f \brackets{\log\parens{\frac{f(X_{i})}{g^{\prime}(X_{i})} \Big/ \frac{f(X_{i})}{g(X_{i})}}}
    \tag{which is well-defined by \Cref{nrmk:log-likelihood-well-defined}.}                                         \\
     & = \EE_f \brackets{\log \frac{f(X_{i})}{g^{\prime}(X_{i})}} - \EE_f \brackets{\log \frac{f(X_{i})}{g(X_{i})}}
    \nonumber                                                                                                       \\
     & = \kldiva{f}{g^{\prime}} - \kldiva{f}{g}.
    \label{neqn:log-likelihood-bernstein-03}
\end{align}
Where the last line follows by definition using \eqref{neqn:kldiv-density-defn}.
We then have
\begin{align}
    \PP( \psi(g,g^{\prime},X) > 0)
     & = \PP\parens{\frac{1}{n}\sum_{i = 1}^{n} \log\frac{g(X_{i})}{g^{\prime}(X_{i})} > 0}
    \tag{using \eqref{neqn:log-likelihood-g-g-prime}}                                                                                                                     \\
     & = \PP\parens{\frac{1}{n}\sum_{i = 1}^{n} \log\frac{g(X_{i})}{g^{\prime}(X_{i})}
        - \EE_f \log \frac{g(X_{1})}{g^{\prime}(X_{1})}
        > \EE_f \log \frac{g^{\prime}(X_{1})}{g(X_{1})}}
    \nonumber                                                                                                                                                             \\
     & = \PP\parens{\frac{1}{n}\sum_{i} \log \frac{g(X_{i})}{g^{\prime}(X_{i})} - \EE_f \log \frac{g(X_{1})}{g^{\prime}(X_{1})} > \kldiva{f}{g} - \kldiva{f}{g^{\prime}}}
    \tag{using \eqref{neqn:log-likelihood-bernstein-03}}                                                                                                                  \\
     & \leq
    \exp \parens{-\frac{n^2 t^{2}}{2\braces{\sum_{i=1}^{n} \EE\brackets{Y_{i}^{2}}+\frac{1}{3} n \kappa t}}}
    \nonumber                                                                                                                                                             \\
     & =
    \exp \parens{-\frac{n^2 t^{2}}{2\braces{n \EE\brackets{Y_{1}^{2}}+\frac{1}{3} n \kappa t}}}
    \tag{since $Y_{i}$ are \iid}                                                                                                                                          \\
     & =
    \exp \parens{-\frac{n t^{2}}{2\braces{ \EE\brackets{Y_{1}^{2}}+\frac{1}{3} \kappa t}}}
    \label{neqn:log-likelihood-bernstein-03b}
\end{align}
where $\kappa \defined 2 \log \beta/\alpha$, $t \defined \kldiva{f}{g} -
\kldiva{f}{g^{\prime}}$, and $Y_{i} \defined \log
\frac{g(X_{i})}{g^{\prime}(X_{i})} - \EE_f \log
\frac{g(X_{1})}{g^{\prime}(X_{1})}$. This follows by the boundedness of $\log
g(X_{i})/g^{\prime}(X_{i})$, and then by applying Bernstein's inequality,
provided that $t > 0$. In order to check this final positivity condition, we
first note that there exists a $C > 0$ such that $\normb{g - g^{\prime}}_{2}
\geq C\delta$, and $\normb{g^{\prime}-f}_{2} \leq \delta$ both hold. We then
have
\begin{equation}\label{neqn:log-likelihood-bernstein-04}
    \normb{f - g}_{2}
    \geq (C - 1) \delta.
\end{equation}
To see this we observe that by assumption, and the triangle inequality
respectively that $C \delta \leq \normb{g - g^{\prime}} \leq \normb{g - f} +
\normb{f - g^{\prime}}$. Then using $\normb{f - g^{\prime}} \leq \delta$ by
assumption and re-arranging, we obtain \eqref{neqn:log-likelihood-bernstein-04}
as required. As a result we obtain the following two inequalities
\begin{align}
    \sqrt{\kldiva{f}{g}}
     & \geq \sqrt{c(\alpha,\beta)} \|f-g\|_2
    \geq
    \sqrt{c(\alpha,\beta)} (C-1) \delta
    \label{neqn:log-likelihood-bernstein-05} \\
    \sqrt{\kldiva{f}{g^{\prime}}}
     & \leq
    \sqrt{1 / \alpha}\|f-g^{\prime}\|_2 \leq \sqrt{1 / \alpha}\delta
    \label{neqn:log-likelihood-bernstein-06},
\end{align}
where $C > 0$ is defined to be a constant satisfying $c(\alpha,\beta) (C-1)^{2}
    > 1 / \alpha$, \ie,
\begin{equation}\label{neqn:crit-cond-bernstein}
    C > 1 + \sqrt{1 / (\alpha c(\alpha,\beta))}.
\end{equation}
Under the condition specified by \eqref{neqn:crit-cond-bernstein}, and by
squaring and subtracting \eqref{neqn:log-likelihood-bernstein-06} from
\eqref{neqn:log-likelihood-bernstein-05}, we obtain
\begin{equation}\label{neqn:log-likelihood-bernstein-07}
    t
    \defined
    \kldiva{f}{g} - \kldiva{f}{g^{\prime}}
    \geq
    (c(\alpha,\beta)(C-1)^{2} - 1 / \alpha)\delta^2
    > 0
\end{equation}
Now we show that $\EE_{f}\parens{Y_{1}^{2}} \lesssim \kldiva{f}{g} +
    \kldiva{f}{g^{\prime}}$. To see this
\begin{align}
    \EE_{f}\parens{Y_{1}^{2}}
     & \leq
    \EE_f \brackets{\parens{\log \frac{g(X_{1})}{g^{\prime}(X_{1})}}^{2}}
    \nonumber                                                               \\
     & =
    \EE_f \brackets{\log\parens{\frac{f(X_{1})}{g^{\prime}(X_{1})} \Big/ \frac{f(X_{1})}{g(X_{1})}}^{2}}
    \tag{which is well-defined by \Cref{nrmk:log-likelihood-well-defined}.} \\
     & =
    \EE_f \brackets{\parens{\log \frac{f(X_{1})}{g^{\prime}(X_{1})} -
            \log \frac{f(X_{1})}{g(X_{1})}}^{2}}
    \nonumber                                                               \\
     & \leq
    2\ubrace{\EE_f \brackets{\parens{\log \frac{f(X_{1})}{g(X_{1})}}^{2}}}{\defines A}{} +
    2\ubrace{\EE_f \brackets{\parens{\log \frac{f(X_{1})}{g^{\prime}(X_{1})}}^{2}}}{\defines B}{}
    \tag{using $(a - b)^{2} \leq 2(a^{2} + b^{2})$, for $a, b \geq 0$.}
\end{align}
We now bound the $A$ term above, with $B$ handled similarly. We observe that:
\begin{align}
    A
     & \defined
    \EE_f \brackets{\parens{\log \frac{f(X_{1})}{g(X_{1})}}^{2}}
    \tag{by definition} \\
     & =
    \int f \parens{\log \frac{f}{g}}^{2} \dlett{\mu}
    \nonumber           \\
     & =
    \int_{f \leq g} f \parens{\log \frac{g}{f}}^{2} \dlett{\mu} +
    \int_{g < f} f \parens{\log \frac{f}{g}}^{2} \dlett{\mu}.
    \label{neqn:bernstein-bdd-second-moment-01}
\end{align}
Now using $\log x \leq x - 1$, for each $x \in \RR_{>0}$, we have that
\begin{equation}\label{neqn:bernstein-bdd-second-moment-02}
    \parens{\log \frac{g}{f}}^{2}
    \leq
    \parens{\frac{g - f}{f}}^{2}
    \text{ and }
    \parens{\log \frac{f}{g}}^{2}
    \leq
    \parens{\frac{f - g}{g}}^{2},
\end{equation}
which hold for $f \leq g$ (\ie, $\frac{g}{f} \geq 1$), and $g < f$ (\ie,
$\frac{f}{g} > 1$), respectively. Now we have:
\begin{align}
    A
     & \leq
    \int_{f \leq g} \frac{(g - f)^{2}}{f} \dlett{\mu} +
    \int_{g < f} \frac{(f - g)^{2} f}{g^{2}} \dlett{\mu}
    \tag{using \eqref{neqn:bernstein-bdd-second-moment-01} and \eqref{neqn:bernstein-bdd-second-moment-02}.} \\
     & \leq
    (1 / \alpha) \int_{f \leq g} (g - f)^{2} \dlett{\mu} +
    (\beta / \alpha^{2}) \int_{g < f} (f - g)^{2} \dlett{\mu}
    \tag{since $0 < \alpha < f, g \leq \beta$}                                                               \\
     & \leq
    (\beta / \alpha^{2}) \normb{f - g}_{2}^{2}
    \tag{since $\beta / \alpha^{2} \geq 1 / \alpha$.}                                                        \\
     & \leq
    K(\alpha, \beta) \kldiva{f}{g},
    \label{neqn:bernstein-bdd-second-moment-03}
\end{align}
where $K(\alpha, \beta) \defined \beta / (\alpha^{2} c(\alpha, \beta))$, where
$c(\alpha, \beta)$ is as defined in \Cref{nlem:equiv-kl-euc-metric}. By a
similar argument, we also have that
\begin{equation}\label{neqn:bernstein-bdd-second-moment-04}
    B
    \leq
    K(\alpha, \beta) \kldiva{f}{g^{\prime}}.
\end{equation}
Let $z \defined \kldiva{f}{g} + \kldiva{f}{g^{\prime}}$. Then using
\eqref{neqn:bernstein-bdd-second-moment-03} and
\eqref{neqn:bernstein-bdd-second-moment-04}, we obtain
\begin{equation}\label{neqn:bernstein-bdd-second-moment-05}
    \EE_{f}\parens{Y_{1}^{2}}
    \leq 2 K(\alpha, \beta) [\kldiva{f}{g} + \kldiva{f}{g^{\prime}}]
    \defines 2 z K(\alpha, \beta)
\end{equation}
Now we use the basic inequality $a + b \leq \parens{\sqrt{a} + \sqrt{b}}^{2}
    \leq 2 (a + b)$, to obtain
\begin{equation}
    z
    \leq
    \parens{\sqrt{\kldiva{f}{g}} + \sqrt{\kldiva{f}{g^{\prime}}}}^{2}
    \leq
    2 z.
\end{equation}
Now, $t^{2} \defined \parens{\kldiva{f}{g} - \kldiva{f}{g^{\prime}}}^{2} =
    \parens{\sqrt{\kldiva{f}{g}} - \sqrt{\kldiva{f}{g^{\prime}}}}^{2}
    \parens{\sqrt{\kldiva{f}{g}} + \sqrt{\kldiva{f}{g^{\prime}}}}^{2}$, we have:
\begin{equation}
    \parens{\sqrt{\kldiva{f}{g}} - \sqrt{\kldiva{f}{g^{\prime}}}}^{2}
    z
    \leq
    t^{2}
    \leq
    2 \parens{\sqrt{\kldiva{f}{g}} - \sqrt{\kldiva{f}{g^{\prime}}}}^{2}
    z.
\end{equation}
We then conclude using
\eqref{neqn:log-likelihood-bernstein-03b},\eqref{neqn:log-likelihood-bernstein-07},
that
\begin{align}
    \PP( \psi(g,g^{\prime},X) > 0)
     & \leq
    \exp \parens{-\frac{n t^{2}}{2\braces{ \EE\brackets{Y_{1}^{2}}+\frac{1}{3} \kappa t}}}
    \tag{per \eqref{neqn:log-likelihood-bernstein-03b}}                                                   \\
     & \leq
    \exp \parens{-\frac{n \parens{\sqrt{\kldiva{f}{g}} - \sqrt{\kldiva{f}{g^{\prime}}}}^{2}
            z}{2\braces{ 2 z K(\alpha, \beta) +\frac{1}{3} \kappa z}}}
    \tag{since $t \leq z$ and \eqref{neqn:bernstein-bdd-second-moment-05}}                                \\
     & =
    \exp \parens{-\frac{n \parens{\sqrt{\kldiva{f}{g}} - \sqrt{\kldiva{f}{g^{\prime}}}}^{2}
        }{2\braces{ 2 K(\alpha, \beta) + \frac{1}{3} \kappa}}}
    \nonumber                                                                                             \\
     & \leq \exp \parens{-\frac{n \parens{\sqrt{c(\alpha,\beta)} (C-1)  -  \sqrt{1 / \alpha}}^{2}\delta^2
        }{2\braces{ 2 K(\alpha, \beta) +\frac{1}{3} \kappa}}}
    \tag{by subtracting
        \eqref{neqn:log-likelihood-bernstein-06} from
    \eqref{neqn:log-likelihood-bernstein-05}}                                                             \\
     & \defines \exp\parens{-n L(\alpha, \beta, C) \delta^{2}},
    \nonumber
\end{align}
whenever condition \eqref{neqn:crit-cond-bernstein} holds, and $L(\alpha, \beta,
    C) \defined \frac{ \parens{\sqrt{c(\alpha,\beta)} (C-1) - \sqrt{1 /
    \alpha}}^{2} }{2\braces{ 2 K(\alpha, \beta) +\frac{2}{3} \log \beta /
    \alpha}}$. Now, taking the supremum over all $g, g^{\prime} \colon \normb{g
    - g^{\prime}}_{2} \geq C\delta, \normb{g^{\prime}-f}_{2} \leq \delta$, the
    required result follows.
\eprfof

\subsection{Proof of
    \texorpdfstring{\Cref{nlem:critical-log-likelihood-concentration}}{\autoref{nlem:critical-log-likelihood-concentration}}}\label{app:nlem:critical-log-likelihood-concentration}
    Recall that \Cref{nlem:critical-log-likelihood-concentration} is concerning
    a packing set. Suppose we have a maximal packing set of $\cF^{\prime}
    \subset \cF$, \ie, $\theset{g_1, \ldots, g_m} \subset \cF^{\prime} \subset
    \cF$ such that $\normb{g_i- g_{j}}_{2} > \delta$ for all $i \neq j$, and it
    is known that $f \in \cF^{\prime}$. We then obtain a key concentration
    result as per \Cref{nlem:critical-log-likelihood-concentration}.
\criticalloglikelihoodconcentration*
\bprfof{\Cref{nlem:critical-log-likelihood-concentration}} We first define the
intermediate \emph{thresholding} random variables
\begin{align*}
    T_{k} \defined
    \begin{cases}
        \max_{j \in [m]} \normb{g_{j}- g_{k}}_{2}
          & , \mbox{s.t. }
        \sum_{i = 1}^{n} \log g_{j}(X_{i}) \geq  \sum_{i = 1}^{n} \log g_{k}(X_{i}), \normb{g_{j}- g_{k}}_{2} > C\delta \\
        0 & , \mbox{otherwise},
    \end{cases}
\end{align*}
for each $k \in [m]$. Without loss of generality suppose that $\normb{g_{k} -
        f}_{2} \leq \delta$. Next
\begin{align*}
    \PP(\|g_{j^*}-f\|_2 > (C + 1)\delta)
     & \leq
    \mathbb{P}(j^* \in \{j: \normb{g_{j} - g_{k}}_{2} > C\delta\})
    \\
     & \leq
    \PP(T_{k} > 0).
\end{align*}
On the other hand
\begin{align*}
    \PP(T_{k} > 0)
     & =
    \PP\parens{\exists j \in [m] \colon
        \sum_{i = 1}^{n} \log g_{j}(X_{i})
        \geq
        \sum_{i = 1}^{n} \log g_{k}(X_{i}),
        \normb{g_{j}- g_{k}}_{2}
    > C \delta}  \\
     & =
    \PP\parens{\bigcup_{j = 1}^{m}
        \braces{\sum_{i = 1}^{n} \log g_{j}(X_{i})
            \geq
            \sum_{i = 1}^{n} \log g_{k}(X_{i}),
            \normb{g_{j}- g_{k}}_{2}
    > C \delta}} \\
     & \leq
    m \exp\parens{-n L(\alpha, \beta, C) \delta^{2}},
    \tag{using union bound and \Cref{nlem:log-likelihood-bernstein}.}
\end{align*}
where $C$ is assumed to satisfy \eqref{neqn:log-likelihood-bernstein-02a}, and
$L(\alpha, \beta, C)$ is defined as per
\eqref{neqn:log-likelihood-bernstein-02b}.
\eprfof

\subsection{Proof of
    \texorpdfstring{\Cref{nlem:local-metric-ent-monotone}}{\autoref{nlem:local-metric-ent-monotone}}}\label{app:nlem:local-metric-ent-monotone}

\localmetricentmonotone*
\bprfof{\Cref{nlem:local-metric-ent-monotone}} It suffices to show that if $g_1,
    \ldots, g_m \in \cF \cap B_{2}(\theta, \varepsilon)$ is a maximal packing
    set at a distance $\varepsilon/c$, then we can pack $ B_{2}(\theta,
    \varepsilon') \cap \cF$ at a distance $\varepsilon'/c$ with at least $m$
    points where $\varepsilon' < \varepsilon$. Consider the points $\theta (1-
    \varepsilon'/\varepsilon) + \varepsilon'/\varepsilon g_{j}$. These points
    clearly are densities since $\theta, g_j \in \cF$. We will show that these
    points are an $\varepsilon'/c$ packing of  $B_{2}(\theta, \varepsilon') \cap
    \cF$. First let us convince ourselves that the points belong to the set. We
    have
\begin{align*}
    \|\theta (1- \varepsilon'/\varepsilon) + \varepsilon'/\varepsilon g_{j} - \theta\|_2 = \varepsilon'/\varepsilon \|g_{j}- \theta\|_2 \leq\varepsilon',
\end{align*}
and using the fact that $\cF$ is convex (by assumption) grants the conclusion.
Next
\begin{align*}
    \|\theta (1-\varepsilon'/\varepsilon) + \varepsilon'/\varepsilon g_{j}- \theta (1-
    \varepsilon'/\varepsilon) - \varepsilon'/\varepsilon g_k\|_2 =
    \varepsilon'/\varepsilon \|g_{j}- g_k\|_2 > \varepsilon'/c,
\end{align*}
which completes the proof.
\eprfof

\subsection{Proof of
    \texorpdfstring{\Cref{nprop:nu-star-estimator-measurability}}{\autoref{nprop:nu-star-estimator-measurability}}}\label{app:nprop:nu-star-estimator-measurability}
\estimatormeasurability*
\bprfof{\Cref{nprop:nu-star-estimator-measurability}} Recall our multistage
sieve estimator $\nu^{*}(\combdatavec) \defined \Upsilon_{\overline{J}}$, where
$\combdatavec \defined (X_{1}, \ldots, X_{n})^{\top}$ is a fixed data sample.
Here $\Upsilon_{\overline{J}}$ denotes the last term of the \emph{finite}
sequence $\Upsilon \defined \theseqb{\Upsilon_{k}}{k = 1}{\overline{J}}$ as
described in \Cref{sec:upper-bound}.

In order to show the measurability of $\nu^{*}(\combdatavec)$ we need to
formalize our setting. We note that our estimator $\nu^{*} \colon B^{n} \to
\cF$, is more precisely a map from the measurable space $(B^{n}, \sigma(B^{n}))$
to the measurable space $(\cF, \sigma(\cF))$. Here $\sigma(B^{n})$ and
$\sigma(\cF)$ denote the Borel $\sigma$-field with respect to the Euclidean and
$L_{2}$-metric topologies on $B^{n}$ and $\cF$, respectively.

Our proof strategy will be to proceed by induction on $k \in [\overline{J}]$
over the sequence $\Upsilon$. We will show that each
$k^{\textnormal{th}}$-indexed map in $\Upsilon$, \ie, $\Upsilon_{k}$, is Borel
measureable, which in turn will imply the measureability of the
$\nu^{*}(\combdatavec)$. Following our (maximal) packing set construction as
described in \Cref{sec:upper-bound} and
\Cref{fig:packing-set-tree-construction}, we need to consider the case where the
traversal down the tree is not necessarily unique at each level, \ie, there may
be collisions (ties) in the packing set children nodes, where the likelihood is
equal. We do always ensure a unique path down the maximal packing set tree, by
selecting the smallest alphanumerically indexed children node at each level.
However, our measurability proof must account for this selection rule
explicitly.

In order to proceed by induction, we consider the base case for $k = 1$, \ie,
$\Upsilon_{1} \in \cF$.  Importantly, we note that $\Upsilon_{1}$ is chosen
arbitrarily from $\cF$ independently of the data samples, $\combdatavec$. Let $A
\in \sigma(\cF)$ be any Borel set. Since all samples $\combdatavec \in B^{n}$
are mapped to $\Upsilon_{1}$ in our setting, then $\Upsilon_{1}^{-1}[A] = B^{n}$
if $\Upsilon_{1} \in A \in \sigma(\cF)$, or $\Upsilon_{1}^{-1}[A] =
\varnothing$, otherwise. In either case we have $\varnothing, B^{n} \in
\sigma(B^{n})$, which shows that $\Upsilon_{1}$ is Borel measurable. Now
consider the event $\theset{\Upsilon_{2} = m_{s}} \defined \thesetb{(X_{1},
\ldots, X_{n})^{\top} \in B^{n}}{\Upsilon_{2}(X_{1}, \ldots, X_{n}) = m_{s}}
\subset B^{n}$, for some index $s \in \NN$. Then we have
\begin{align}
    \theset{\Upsilon_{2} = m_{s}}
     & \defined
    \thesetb{(X_{1}, \ldots, X_{n})^{\top} \in B^{n}}{\Upsilon_{2}(X_{1}, \ldots, X_{n}) = m_{s}}
    \label{neqn:upsilon2-measurability-01}
    \\
     & =
    \bigcap_{g \in P_{\Upsilon_{1}}}
    \thesetb{(X_{1}, \ldots, X_{n})^{\top} \in B^{n}}{\sum_{i = 1}^{n}
        \log (m_{s}(X_{i})) \geq \sum_{i = 1}^{n} \log (g(X_{i}))} \; \bigcap
    \nonumber                  \\
     & \bigcap_{j = 1}^{s - 1}
    \thesetb{(X_{1}, \ldots, X_{n})^{\top} \in B^{n}}{\sum_{i = 1}^{n}
        \log (m_{s}(X_{i})) > \sum_{i = 1}^{n} \log (m_{j}(X_{i}))}.
    \label{neqn:upsilon2-measurability-02}
\end{align}
In \eqref{neqn:upsilon2-measurability-02}, we observe that $\theset{\Upsilon_{2}
        = m_{s}} \subset B^{n}$ is represented as the intersection of 2 separate
        (finite) set intersections. Note that the second intersection set
        \emph{explicitly} accounts for our alphanumerical index selection rule
        in the children densities of $P_{\Upsilon_{1}}$. Consider the first
        finite intersection term. Here, each $g \in P_{\Upsilon_{1}} \subset
        \cF$ are Borel measurable by \eqref{neqn:density-class-F-alpha-beta}. We
        note that the $\log$ and the addition (\ie, ``$+_{\RR}$'') functions are
        both continuous and measurable, and therefore, so is their composition.
        Thus the resulting \emph{finite} sum, $\sum_{i = 1}^{n} \log f(X_{i})$,
        is a measurable function, for any density $f \in \cF$ (which is always
        measurable). As such the $\Upsilon_{2}$ is measurable since all these
        inequalities give rise to measurable sets and when one intersects them
        (they are finitely many) one obtains another measurable set. Once again,
        let $A \in \sigma(\cF)$ be any Borel set. Then such an $A$ contains
        either no such densities $m_{s}$, or at most finitely many (since the
        number of children of our maximal packing set tree is always finite). If
        no such $m_{s} \in A$, then $\Upsilon_{2}^{-1}[A] = \varnothing \in
        \sigma(B^{n})$. Thus $\Upsilon_{2}$ is indeed Borel measurable in this
        case. In the case where there exist finitely many such $m_{s} \in A$, it
        follows that
\begin{equation}\label{neqn:upsilon2-measurability-03}
    \Upsilon_{2}^{-1}[A]
    =
    \bigcup_{\thesetb{s}{m_{s} \in A}} \theset{\Upsilon_{2} = m_{s}}
    \defines
    \bigcup_{\thesetb{s}{m_{s} \in A}} \thesetb{(X_{1}, \ldots, X_{n})^{\top} \in B^{n}}{\Upsilon_{2}(X_{1}, \ldots, X_{n}) = m_{s}}
\end{equation}
In \eqref{neqn:upsilon2-measurability-03} we note that $\Upsilon_{2}^{-1}[A]$
represents a finite union of Borel measurable sets as per
\eqref{neqn:upsilon2-measurability-02}, which is again Borel measurable. That
is, we have shown that $\Upsilon_{2}^{-1}[A] \in \sigma(B^{n})$, which indeed
implies the Borel measurability of $\Upsilon_{2}$, as required.

Similarly, consider the event $\theset{\Upsilon_{3} \defines m_{s, t}} \subset
    B^{n}$, for some $t \in \NN$ and $s \in \NN$ taken as per
    \eqref{neqn:upsilon2-measurability-01}. Here, the indexed density $m_{s, t}$
    signifies that $\Upsilon_{3}$ is derived from the children of the packing
    set of $\Upsilon_{2} \defines m_{s}$, as denoted by $P_{m_{s}}$ in our work.
    Once again we can write this $\Upsilon_{3}$ as
\begin{align}
    \theset{\Upsilon_{3} = m_{s, t}}
     & \defined
    \thesetb{(X_{1}, \ldots, X_{n})^{\top} \in B^{n}}{\Upsilon_{3}(X_{1}, \ldots, X_{n}) = m_{s, t}}
    \nonumber                                     \\
     & =
    \bigcap_{g \in P_{m_{s}}}
    \thesetb{(X_{1}, \ldots, X_{n})^{\top} \in B^{n}}{\sum_{i = 1}^{n}
        \log (m_{s, t}(X_{i})) \geq \sum_{i = 1}^{n} \log (g(X_{i}))} \; \bigcap
    \nonumber                                     \\
     & \bigcap_{j = 1}^{t - 1}
    \thesetb{(X_{1}, \ldots, X_{n})^{\top} \in B^{n}}{\sum_{i = 1}^{n}
        \log (m_{s, t}(X_{i})) > \sum_{i = 1}^{n} \log (m_{s, j}(X_{i}))}
    \; \bigcap \nonumber                          \\
     & \bigcap \;  \theset{\Upsilon_{2} = m_{s}}.
\end{align}
By a similar argument to the measurability of $\Upsilon_{2}$ it follows that
$\Upsilon_{3}$ is also measurable. As such, given the recursive construction of
the finite sequence $\Upsilon \defined \theseqb{\Upsilon_{k}}{k =
1}{\overline{J}}$ via our maximal packing set tree traversal, this pattern
inductively repeats for each $k \in \theset{4, \ldots, \overline{J}}$. Since
$\nu^{*}(\combdatavec) \defined \Upsilon_{\overline{J}}$, this implies the
measurability of $\nu^{*}(\combdatavec)$, as required.
\eprfof

\begin{restatable}{nrmk}{countablymeasureable}\label{nrmk:countably-measureable}
We note that the arguments in the proof above hold, even if the cardinality of
the set of children densities at any iteration were at most countable (not
\emph{just} finite). That is, \eqref{neqn:upsilon2-measurability-02} would still
return a measureable set even if $\absa{P_{\Upsilon_{1}}} = \infty$, since Borel
measurability is preserved over countable intersections and unions. The packing
sets in our construction are necessarily at most countable, since all of the
subsets of $\cF$ we consider are separable in the $L_{2}$-metric (\ie, contain a
countably dense subset). This follows from \Cref{nlem:separability-of-class-F}.
\end{restatable}

\subsection{Proof of
    \texorpdfstring{\Cref{nthm:upper-bound-rate-finite-iterations}}{\autoref{nthm:upper-bound-rate-finite-iterations}}}\label{app:nthm:upper-bound-rate}

We begin with a useful result, which will enable us to construct upper bounds
for estimator $\nu^{*}(\combdatavec)$.
\begin{restatable}{nlem}{upsiloniscauchy}\label{nlem:upsilon-is-cauchy-sequence}
    The finite sequence $\Upsilon \defined \theseqb{\Upsilon_{k}}{k =
    1}{\overline{J}}$, as defined in the construction of our estimator
    $\nu^{*}(\combdatavec)$, satisfies
    \begin{equation}
        \normb{\Upsilon_{J} - \Upsilon_{J^{\prime}}}_{2}
        \leq
        \frac{d}{2^{J^{\prime} - 2}},
    \end{equation}
    for each pair of positive integers $J^{\prime} < J$.
\end{restatable}
\bprfof{\Cref{nlem:upsilon-is-cauchy-sequence}} Let $\Upsilon_{J^{\prime}},
    \Upsilon_{J} \in \Upsilon$, for any positive integers $J > J^{\prime} \geq
    1$. We then have
\begin{equation}\label{neqn:upsilon-is-cauchy-sequence-01}
    \normb{\Upsilon_{J} - \Upsilon_{J^{\prime}}}_{2}
    \leq
    \sum_{i = J^{\prime}}^{J - 1}
    \normb{\Upsilon_{i + 1} - \Upsilon_{i}}_{2}
    \leq
    \sum_{i = J^{\prime}}^{J - 1} \frac{d}{2^{i-1}}
    \leq
    \frac{d}{2^{J^{\prime} - 2}}.
\end{equation}
As required.
\eprfof
\begin{restatable}[Telescoping sum of conditional
        probabilities]{nlem}{telescopingsumscondprobs}\label{nlem:telescoping-sum-cond-prob}
        Let $n \geq 2$ be a fixed integer, and $\theset{A_{1}, A_{2}, \ldots,
        A_{n}}$ denote events on a common probability space, with
        $\PP(\stcomp{A}_{j}) > 0$ for each $j \geq 1$. We then have
    \begin{equation}\label{neqn:telescoping-sum-cond-prob-01}
        \PP(A_{n})
        \leq
        \sum_{j = n}^{2}\PP(A_{j} \mid \stcomp{A}_{j -1}) + \PP(A_{1}).
    \end{equation}
\end{restatable}
\bprfof{\Cref{nlem:telescoping-sum-cond-prob}} We will prove this by induction
on $n \geq 2$. We check the induction base case for $n = 2$. We first observe
that
\begin{equation}\label{neqn:telescoping-sum-cond-prob-02}
    A_{2}
    \subseteq
    A_{1} \cup A_{2}
    =
    (A_{2} \cap \stcomp{A}_{1})
    \sqcup
    A_{1},
\end{equation}
where the latter set is a \emph{disjoint} union. It then follows that
\begin{align}
    \PP(A_{2})
     & \leq
    \PP\parens{A_{2} \cap \stcomp{A}_{1}} +
    \PP(A_{1})
    \tag{by monotonicity of $\PP$ applied to  \eqref{neqn:telescoping-sum-cond-prob-02}} \\
     & \leq
    \frac{\PP\parens{A_{2} \cap \stcomp{A}_{1}}}{\PP(\stcomp{A}_{1})} +
    \PP(A_{1})
    \tag{since $\PP(\stcomp{A}_{1}) \in (0, 1]$, by assumption}                          \\
     & \defines
    \PP(A_{2} \mid \stcomp{A}_{1}) + \PP(A_{1}),
    \nonumber
\end{align}
which proves the base case for $n = 2$. Now, by induction assume the result is
true for each integer $n = k > 2$. We then have for $n = k + 1$ that:
\begin{align}
    \PP(A_{k + 1})
     & \leq
    \PP(A_{k + 1} \mid \stcomp{A}_{k}) + \PP(A_{k})
    \tag{using induction base case}  \\
     & \leq
    \PP(A_{k + 1} \mid \stcomp{A}_{k}) +
    \sum_{j = k}^{2}\PP(A_{j} \mid \stcomp{A}_{j -1}) + \PP(A_{1})
    \tag{using induction hypothesis} \\
     & =
    \sum_{j = k + 1}^{2}\PP(A_{j} \mid \stcomp{A}_{j -1}) + \PP(A_{1}),
\end{align}
as required. So the result is true for $n = k + 1$, and thus by induction holds
for each integer $n \geq 2$.
\eprfof

\upperboundratefiniteiters*
\bprfof{\Cref{nthm:upper-bound-rate-finite-iterations}} Combining the results of
Lemma \ref{nlem:critical-log-likelihood-concentration} (with $c \defined 2(C+1)$
where $c$ is the constant from the definition of local packing entropy) and
Lemma \ref{nlem:local-metric-ent-monotone} we conclude that for each $j \in
\theset{2, \ldots, J}$ we have
\begin{align}
     & \mathbb{P}\parens{\normb{f - \Upsilon_{j}}_{2}
        >
        \frac{d}{2^{j-1}}
        \,\middle|\,
        \normb{f - \Upsilon_{j-1}}_{2}
        \leq
        \frac{d}{2^{j-2}},
        \Upsilon_{j - 1}}
    \nonumber
    \\
     & \leq
    \absb{P_{\Upsilon_{j-1}}}
    \exp\parens{-\frac{n L(\alpha, \beta, C) d^2}{2^{2(j-1)}(C+1)^2}}
    \label{neqn:telescoping-bound-01}                 \\
     & \leq
    \mlocc{\cF}{\frac{d}{2^{J - 2}}}{c}
    \exp\parens{-\frac{n L(\alpha, \beta, C) d^2}{2^{2(j-1)}(C+1)^2}}
    \label{neqn:telescoping-bound-02}
\end{align}
where $P_{\Upsilon_{j}}$ are the maximal packing sets described in the
construction of $\nu^{*}(\combdatavec)$. Crucially, we observe that the $\rhs$
of \eqref{neqn:telescoping-bound-02} does not depend on the conditioned random
variables, \ie, $\Upsilon_{j - 1}$, for each $j \in \theset{2, \ldots, J}$ hence
we can drop $\Upsilon_{j - 1}$ from the conditioning. Now let denote $A_{j}
\defined \theset{\normb{f - \Upsilon_{j}}_{2} > \frac{d}{2^{j-1}}}$, for each
integer $j \geq 1$. Then we can proceed by working with the unconditional events
$A_{j}$ in \eqref{neqn:telescoping-bound-02}.

Moreover, we then have that $\stcomp{A}_{j - 1} \defined \theset{\normb{f -
            \Upsilon_{j - 1}}_{2}
        \leq
        \frac{d}{2^{j-2}}}$ for each integer $j \geq 2$. In particular
$\PP(\stcomp{A}_{1}) = \theset{\normb{f - \Upsilon_{1}}_{2} \leq d} = 1$, since
$f, \Upsilon_{1} \in \cF$, so indeed $\normb{f - \Upsilon_{1}}_{2} \leq
\diam_{2}{(\cF)} \defines d$ almost surely. By aligning our notation directly
with \Cref{nlem:telescoping-sum-cond-prob}, we can apply the telescoping bound
to $\PP(A_{j})$ as follows
\begin{align}
    \PP(A_{J})
     & \defined
    \mathbb{P}\parens{\normb{f - \Upsilon_{J}}_{2}
        >
        \frac{d}{2^{J-1}}}
    \tag{by definition}                         \\
     & \leq
    \mlocc{\cF}{\frac{d}{2^{J - 2}}}{c}
    \sum_{j = 1}^{J-1}\exp\parens{-\frac{n L(\alpha, \beta, C) d^2}{2^{2j}(C+1)^2}}
    \tag{per \eqref{neqn:telescoping-bound-02}} \\
     & \leq
    \mlocc{\cF}{\frac{d}{2^{J - 2}}}{c}
    a (1 + a^{4-1} + a^{16-1} + \ldots)\mathbbm{1}(J > 1)
    \label{neqn:telescoping-bound-03}           \\
     & \leq
    \mlocc{\cF}{\frac{d}{2^{J - 2}}}{c}
    a (1 + a + a^{2} + \ldots)\mathbbm{1}(J > 1)
    \nonumber                                   \\
     & \leq
    \mlocc{\cF}{\frac{d}{2^{J - 2}}}{c}
    \frac{a}{1-a} \mathbbm{1}(J > 1),
    \label{neqn:telescoping-bound-04}
\end{align}
where for brevity in \eqref{neqn:telescoping-bound-03} we denote
\begin{align*}
    a
    \defined
    \exp\parens{-\frac{n L(\alpha, \beta, C) d^2}{2^{2(J-1)}(C+1)^2}}.
\end{align*}
Since $C$ is assumed to satisfy \eqref{neqn:log-likelihood-bernstein-02a}, and
$L(\alpha, \beta, C)$ is defined as per
\eqref{neqn:log-likelihood-bernstein-02b}, it follows that $a < 1$. Note here
that the above bound \eqref{neqn:telescoping-bound-04} holds, provided that
$\PP(\stcomp{A}_{j}) > 0$ for $j < J$ as required by
\Cref{nlem:telescoping-sum-cond-prob}. Suppose that the $\rhs$ of
\eqref{neqn:telescoping-bound-04} is strictly smaller than $1$. In that case for
all $j$, $\PP(\stcomp{A}_{j}) > 0$ since bound \eqref{neqn:telescoping-bound-04}
holds inductively for all $\PP(A_j)$ for $j \leq J$. On the other hand, if the
$\rhs$ of \eqref{neqn:telescoping-bound-04} is $\geq 1$ then
\eqref{neqn:telescoping-bound-04} trivially holds. In both cases we conclude
that \eqref{neqn:telescoping-bound-04} holds.

If one sets $\varepsilon_J \defined \frac{\sqrt{L(\alpha, \beta, C)}
        d}{2^{(J-1)}(C+1)}$, we have that if
\begin{align*}
    n \varepsilon_J^2 > 2 \log
    \mlocc{\cF}{\varepsilon_J \frac{2 (C+1)}{\sqrt{L(\alpha, \beta, C)}}}{c} = 2
    \log \mlocc{\cF}{\frac{d}{2^{J - 2}}}{c},
\end{align*} and $a \defined \exp(-n \varepsilon_J^2) <
    1/2 \iff n \varepsilon_J^2 > \log 2$, the above probability in
\eqref{neqn:telescoping-bound-04} will be bounded from above by $2 \exp(-n
\varepsilon_J^2 / 2)$. This condition is implied when
\begin{align}\label{neqn:suff-condition-epsJ}
    n \varepsilon_J^2
    >
    2 \log \mlocc{\cF}{\varepsilon_J
        \frac{2 (C+1)}{\sqrt{L(\alpha, \beta, C)}}}{c}
    \vee \log 2.
\end{align}
We now have
\begin{align}\label{neqn:mu-upislon-ineq}
    \normb{\nu^*_{\overline{J}} - f}_2
    \leq
    \normb{\Upsilon_{\overline{J}} - \Upsilon_J}_2 +
    \normb{\Upsilon_J - f}_2
    \leq
    3 \varepsilon_J \frac{C+1}{\sqrt{L(\alpha, \beta, C)}},
\end{align} with probability at least $1 - 2\exp(-n\varepsilon_J^2 / 2)$ which holds for all $J$ satisfying \eqref{neqn:suff-condition-epsJ} (including $\overline{J}$). Here we want to clarify that the last inequality in \eqref{neqn:mu-upislon-ineq} follows from the fact that $\|\Upsilon_{\overline{J}} - \Upsilon_J\|_2 \leq d/2^{J-2}$, as seen when we verified that $\Upsilon$ forms a Cauchy sequence in \Cref{nlem:upsilon-is-cauchy-sequence} (and since $\overline{J} \geq J$).
Let $J^*$ be selected as the maximum integer $J$ such that
\eqref{neqn:suff-condition-epsJ} holds, or otherwise if such $J$ does not exist
$J^* = 1$, i.e. $J^* \equiv \overline{J}$. Let $\eta =
3\frac{C+1}{\sqrt{L(\alpha, \beta, C)}}$, $\underline{C} = 2$ and $C' = 1 / 2$.
We have established that the following bound holds
\begin{equation*}
    \mathbb{P}(\normb{f - \nu^*_{\overline{J}}}_{2}
    >
    \eta \varepsilon_J)
    \leq
    \underline C \exp(-C'n\varepsilon_J^2) \mathbbm{1}(J > 1)
    \leq
    \underline C \exp(-C'n\varepsilon_J^2) \mathbbm{1}(J^* > 1),
\end{equation*}
for all $1 \leq J \leq J^*$, where this bound also holds in the case when $J^* =
    1$ by exception. Observe that we can extend this bound to all $J \in
    \mathbb{Z}$ and $J \leq J^*$, since for $J < 1$ we have $\eta \varepsilon_J
    \geq 6 d$ and so
\begin{equation*}
    \mathbb{P}(\normb{f - \nu^*_{\overline{J}}}_{2}
    >
    \eta \varepsilon_J)
    \leq
    0
    \leq
    \underline C \exp(-C'n\varepsilon_J^2) \mathbbm{1}(J^* > 1).
\end{equation*}
We conclude that
\begin{equation*}
    \mathbb{P}(\normb{f - \nu^*_{\overline{J}}}_{2}
    >
    \eta \varepsilon_J)
    \leq
    0
    \leq
    \underline C \exp(-C'n\varepsilon_J^2) \mathbbm{1}(J^* > 1),
\end{equation*}
for any $J \leq J^*$. Now for any $\varepsilon_{J-1} > x \geq \varepsilon_{J}$
for $J \leq J^*$ we have that
\begin{align*}
    \mathbb{P}(\normb{f - \nu^*_{\overline{J}}}_{2}
    >
    2\eta x)
     & \leq
    \mathbb{P}(\normb{f - \nu^*_{\overline{J}}}_{2}
    >
    \eta \varepsilon_{J - 1})                                       \\
     & \leq
    \underline C \exp(-C'n\varepsilon_{J-1}^2) \mathbbm{1}(J^* > 1) \\
     & \leq
    \underline C \exp(-C'nx^2)\mathbbm{1}(J^* > 1),
\end{align*}
where the last inequality follows due to the fact that the map $x \mapsto
    \underline C \exp(-C'nx^2)$ is monotonically decreasing for positive reals.
    We will now integrate the tail bound:
\begin{equation}\label{important:prob:bound}
    \mathbb{P}(\normb{f - \nu^*_{\overline{J}}}_{2}
    >
    2 \eta x)
    \leq
    \underline C \exp(-C'nx^2) \mathbbm{1}(J^* > 1),
\end{equation}
which holds true for $x \geq \varepsilon^* \defined \varepsilon_{J^{*}}$, where
    $\varepsilon_J = \frac{\sqrt{L(\alpha, \beta, C)} d}{2^{(J-1)}(C+1)}$,
    always (since even if $J^* = 1$ by exception, this bound is still valid). We
    then have
\begin{align*}
    \EE \normb{f - \nu^*_{\overline{J}}}_{2}^{2}
     & =
    \int_{0}^{\infty} 2 x
    \mathbb{P}(\normb{f - \nu^*_{\overline{J}}}_{2} > x) \dlett{x}                                                            \\
     & \leq
    C''' \varepsilon^{*2} + \int_{2\eta\varepsilon^*}^{\infty} 2 x \underline C \exp(-C''nx^2) \mathbbm{1}(J^* > 1) \dlett{x} \\
     & =
    C''' \varepsilon^{*2} + C^{''''}n^{-1}\exp(-C'''''n\varepsilon^{*2})\mathbbm{1}(J^* > 1).
\end{align*}
Now $n\varepsilon^{*2}$ is bigger than a constant (\ie, $\log 2$) otherwise $J^*
    = 1$. Hence, the above is smaller than $\bar C \varepsilon^{*2}$ for some
    absolute constant $\bar C$.
\eprfof

\subsection{Proof of
    \texorpdfstring{\Cref{nthm:sharp-minimax-rate}}{\autoref{nthm:sharp-minimax-rate}}}\label{app:nthm:sharp-minimax-rate}
\sharpminimaxrate*
\bprfof{\Cref{nthm:sharp-minimax-rate}} First suppose that $\varepsilon^*$
satisfies $n\varepsilon^{*2} > 4 \log 2$. Then for $\delta^* :=
\varepsilon^*/\sqrt{4(1 / \alpha \vee 1)}$ we have $\log
    \mlocc{\cF}{\delta^*}{c}
    \geq
    \log \mlocc{\cF}{\varepsilon^*}{c}
    \geq
    n\varepsilon^{*2}/2 + n \varepsilon^{*2}/2 > 2n \delta^{*2} / \alpha + 2
    \log 2$ and so this implies the sufficient condition for the lower bound per
    \Cref{nlem:minimax-lower-bound}. Let $\eta \defined \frac{\sqrt{2}
    c}{\sqrt{L(\alpha, \beta, c / 2 - 1)}} \wedge 1$. For a constant $C$ such
    that $C \eta > 1$, we have
\begin{align*}
    C^2 n \varepsilon^{*2}
     & \geq
    1/\eta^2
    \log \mlocc{\cF}{C \eta \varepsilon^*}{c}
    \geq
    \log \mlocc{\cF}{C \eta \varepsilon^*}{c} \\
     & \geq
    \log \mlocc{\cF}{\sqrt{2} C \varepsilon^* \frac{c}{\sqrt{L(\alpha, \beta, c / 2 - 1)}}}{c}
\end{align*}
Setting $\delta \defined C \sqrt{2} \varepsilon^*$ we obtain that
\begin{align*}
    n\delta^2
    \geq
    2 \log \mlocc{\cF}{\delta \frac{c}{\sqrt{L(\alpha, \beta, c / 2 - 1)}}}{c}.
\end{align*}
In addition since $C > 1$, $\delta$ satisfies \eqref{upper:bound:suff:cond}
(taking into account that $n \varepsilon^{*2} > 4 \log 2$, which implies $n
\delta^2 \geq 4 \log 2 C^2 > \log 2$). We note that the map $0 < x \mapsto n x^2
- \log \mlocc{\cF}{x \frac{c}{\sqrt{L(\alpha, \beta, c / 2 - 1)}}}{c} \vee \log
2$ is non-decreasing by \Cref{nlem:local-metric-ent-monotone}. Now, with
$\varepsilon_{J^*}$ defined as per
\Cref{nthm:upper-bound-rate-finite-iterations}, this implies that $\delta \geq
\varepsilon_{J^*}/2$. This shows that the rate in this case is of the order
$\varepsilon^{*2}$.

Next, suppose that $\varepsilon^{*}$ defined by $\sup \{\varepsilon: n
    \varepsilon^{2} \leq \log{\mlocc{\cF}{\varepsilon}{c}}\}$ satisfies $n
    \varepsilon^{*2} \leq 4\log 2$. For $2\varepsilon^*$, we have $16 \log 2
    \geq 4\varepsilon^{*2}n \geq \log{\mlocc{\cF}{2 \varepsilon^{*}}{c}}$. If
    $c$ in the definition of local packing is large enough, we could put points
    in the diameter of the ball with radius $2 \varepsilon^*$ such that the
    packing set has more than $\exp(16\log 2)$ many points. But that implies
    that the set $\cF$ is entirely inside a ball of radius $\sqrt{16\log 2}
    n^{-1/2}$ (as $\varepsilon^{*2} \leq (4\log 2) n^{-1} $). To see the latter,
    one can take the midpoint of the line segment connecting the endpoints of a
    diameter of $\cF$ and position a ball of radius $2\varepsilon^*$ there. In
    such a case, for the lower bound, we could pick $\varepsilon$ to be
    proportional to the diameter of the set (with a small proportionality
    constant). That will ensure that $\varepsilon \sqrt{n}$ is upper bounded by
    some constant (as $2\sqrt{(16 \log 2)}n^{-1/2}$ is bigger than the
    diameter), and at the same time $\log{\mlocc{\cF}{\varepsilon}{c}}$ can be
    made bigger than a constant (provided that $c$ in the definition of a local
    packing is large enough) -- by taking $\theta$ (where $\theta$ is the center
    of the localized set $B_2(\theta, \varepsilon) \cap \cF$) to be the midpoint
    of a diameter of the set $\cF$ and then placing equispaced points on the
    diameter. Hence, the diameter of the set is a lower bound (up to constant
    factors) in this case, which is of course always an upper bound too (up to
    constant factors). So we conclude that either for $\varepsilon^{*}$ defined
    by $\sup \thesetc{\varepsilon}{ \varepsilon^{2} n \leq
    \log{\mlocc{\cF}{\varepsilon}{c}}}$ satisfies $\varepsilon^{*2}n > 4\log 2$
    or the lower and upper bounds are of the order of the diameter of the set.
    In summary the rate is given by the $\varepsilon^{*2} \wedge d^2$. This is
    true since in the second case, $4\varepsilon^*$ is bigger than the diameter
    of the set.
\eprfof

\subsection{Proof of
    \texorpdfstring{\Cref{nprop:extend-zero-bounded-densities}}{\autoref{nprop:extend-zero-bounded-densities}}}\label{app:nprop:extend-zero-bounded-densities}

\extendzeroboundeddensities*
\bprfof{\Cref{nprop:extend-zero-bounded-densities}} We argue this as follows.
Let $f_{\alpha} \in \cF$, which is lower bounded by some $\alpha > 0$. Now
consider the following set of $f_{\alpha}$-mixture densities, \ie, $\cF^{\prime}
= \thesetc{(1/2) f_{\alpha} + (1/2) f}{f \in \cF} \subset \cF$. By construction,
all densities in $\cF^{\prime}$ are thus lower bounded by $\alpha/2$, \ie
$\cF^{\prime} \subset \cF_B^{[\alpha/2,\beta]}$. Moreover, $\cF^{\prime}$ forms
a convex density class. Hence, the minimax rate would be given by
$\varepsilon^{2} \wedge \operatorname{diam}_2(\cF^{\prime})^2$ where
$\varepsilon = \sup \thesetc{\varepsilon}{n \varepsilon^{2} \leq
\log{\mlocc{\cF^{\prime}}{\varepsilon}{c}}}$. We can artificially create
variables from the class $\cF^{\prime}$ by randomizing $X_{i}$ as follows
\begin{align*}
    Z_{i} = \begin{cases}
                T_{i} \distiid f_{\alpha} & \mbox{ with probability } $1/2$, \\
                X_{i}                     & \mbox{ with probability } $1/2$.
            \end{cases}
\end{align*}
Then let $\hat f$ be our estimator of $(1/2) f_{\alpha} + (1/2)f$. We know:
\begin{align*}
    \EE_Z \normb{\hat f - ((1/2) f_{\alpha} + (1/2) f)}_{2}^{2} \lesssim \varepsilon^{2} \wedge \operatorname{diam}(\cF^{\prime})^2,
\end{align*}
so that
\begin{align*}
    \EE_{X, T, V} \|(2\hat f-f_{\alpha}) - f\|_{2}^{2} \lesssim 4 \varepsilon^{2} \wedge \operatorname{diam}(\cF^{\prime})^2,
\end{align*}
where $T = (T_1,\ldots, T_n)$ and $V = (V_1,\ldots, V_n)$ are the values of the
coin flips in the definition of $Z_i$. Hence, $\EE_{T,V}2 \hat f-f_{\alpha}$
achieves the same rate for $f$ since by Jensen's inequality
\begin{align*}
    \EE_Y \|\EE_{T, V} (2\hat f-f_{\alpha}) - f\|_{2}^{2} \leq \EE_{Y, T, V} \|(2\hat f-f_{\alpha}) - f\|_{2}^{2} \lesssim 4 \varepsilon^{2} \wedge \operatorname{diam}(\cF^{\prime})^2.
\end{align*}
Moreover, note that since $\hat f \in \cF^{\prime}$ for each of value of $T,V$
we have $\EE_{T,V} 2\hat f-f_{\alpha} \in \cF$. Thus, the upper bound is the
same for the two sets. On the other hand since $\cF^{\prime} \subset \cF$ the
lower bound is also of the same rate. Finally, observe that $
\log{\mlocc{\cF^{\prime}}{\varepsilon}{c}} = \log{\mlocc{\cF}{2
\varepsilon}{c}}$ so that the order of $\varepsilon^* = \sup
\thesetc{\varepsilon}{n \varepsilon^{2} \leq
\log{\mlocc{\cF^{\prime}}{\varepsilon}{c}}}$ is the same as that of the equation
$\varepsilon^* = \sup \thesetc{\varepsilon}{n \varepsilon^{2} \leq
\log{\mlocc{\cF}{\varepsilon}{c}}}$. In addition, it is also clear that
$2\diam_{2}(\cF^{\prime}) = \diam_{2}(\cF)$.
\eprfof

\subsection{Proof of
    \texorpdfstring{\Cref{nthm:adaptive-upper-bound-rate-finite-iterations}}{\autoref{nthm:adaptive-upper-bound-rate-finite-iterations}}}\label{app:nthm:adaptive-upper-bound-rate-finite-iterations}

We first prove the following simple lemma.
\begin{restatable}[]{nlem}{nlemadaptivityentropybound}\label{nlem:adaptivity-entropy-bound}
        Suppose $\nu, \mu \in \cF$ are two densities such that $\normb{\nu -
        \nu}_{2} \leq \delta$. If $\delta \leq \varepsilon$ then have
        $M(\nu,\varepsilon,c) \leq M(\mu,2\varepsilon,2c)$.
\end{restatable}

\bprfof{\Cref{nlem:adaptivity-entropy-bound}} It suffices to show that $B(\nu,
\varepsilon) \subseteq B(\mu, 2\varepsilon)$. For any $x \in B(\nu,
\varepsilon)$ we have $\|x-  \nu\|_2 \leq \varepsilon$, and hence by the
triangle inequality we obtain
\begin{align*}
    \|x-\mu\|_2 \leq \|x - \nu\|_2 +\|\nu - \mu\|_2 \leq \varepsilon + \delta\leq 2\varepsilon,
\end{align*}
which completes the proof.
\eprfof

\adaptiveupperboundratefiniteiters*
\bprfof{\Cref{nthm:adaptive-upper-bound-rate-finite-iterations}}

Combining the results of Lemma \ref{nlem:critical-log-likelihood-concentration}
(with $c \defined 2(C+1)$ where $c$ is the constant from the definition of local
packing entropy) and Lemma \ref{nlem:local-metric-ent-monotone} we conclude that
for each $j \in \theset{2, \ldots, J}$ we have
\begin{align}
     & \mathbb{P}\parens{\normb{f - \Upsilon_{j}}_{2}
        >
        \frac{d}{2^{j-1}}
        \,\middle|\,
        \normb{f - \Upsilon_{j-1}}_{2}
        \leq
        \frac{d}{2^{j-2}},
        \Upsilon_{j - 1}}
    \nonumber
    \\
     & \leq
    \absb{P_{\Upsilon_{j-1}}}
    \exp\parens{-\frac{n L(\alpha, \beta, C) d^2}{2^{2(j-1)}(C+1)^2}}
    \label{neqn:adaptive-telescoping-bound-01}                 \\
         & \leq
    M(f,\frac{d}{2^{j-3}},2c)
    \exp\parens{-\frac{n L(\alpha, \beta, C) d^2}{2^{2(j-1)}(C+1)^2}}
    \label{neqn:adaptive-telescoping-bound-02}                 \\
     & \leq
    M(f,\frac{d}{2^{J-3}},2c)
    \exp\parens{-\frac{n L(\alpha, \beta, C) d^2}{2^{2(j-1)}(C+1)^2}}
    \label{neqn:adaptive-telescoping-bound-03}
\end{align}
where $P_{\Upsilon_{j}}$ are the maximal packing sets described in the
construction of $\nu^{*}(\combdatavec)$. Furthermore, inequality
\eqref{neqn:adaptive-telescoping-bound-02} follows from Lemma
\ref{nlem:adaptivity-entropy-bound}. Crucially, we observe that the $\rhs$ of
\eqref{neqn:adaptive-telescoping-bound-03} does not depend on the conditioned
random variables, \ie, $\Upsilon_{j - 1}$, for each $j \in \theset{2, \ldots,
J}$ hence we can drop $\Upsilon_{j - 1}$ from the conditioning. Now let denote
$A_{j} \defined \theset{\normb{f - \Upsilon_{j}}_{2} > \frac{d}{2^{j-1}}}$, for
each integer $j \geq 1$. Then we can proceed by working with the unconditional
events $A_{j}$ in \eqref{neqn:telescoping-bound-02}.

Moreover, we then have that $\stcomp{A}_{j - 1} \defined \theset{\normb{f -
            \Upsilon_{j - 1}}_{2}
        \leq
        \frac{d}{2^{j-2}}}$ for each integer $j \geq 2$. In particular
$\PP(\stcomp{A}_{1}) = \theset{\normb{f - \Upsilon_{1}}_{2} \leq d} = 1$, since
$f, \Upsilon_{1} \in \cF$, so indeed $\normb{f - \Upsilon_{1}}_{2} \leq
\diam_{2}{(\cF)} \defines d$ almost surely. By aligning our notation directly
with \Cref{nlem:telescoping-sum-cond-prob}, we can apply the telescoping bound
to $\PP(A_{j})$ as follows
\begin{align}
    \PP(A_{J})
     & \defined
    \mathbb{P}\parens{\normb{f - \Upsilon_{J}}_{2}
        >
        \frac{d}{2^{J-1}}}
    \tag{by definition}                         \\
     & \leq
    \madlocc{\cF}{f}{\frac{d}{2^{J-3}}}{2c}
    \sum_{j = 1}^{J-1}\exp\parens{-\frac{n L(\alpha, \beta, C) d^2}{2^{2j}(C+1)^2}}
    \tag{per \eqref{neqn:telescoping-bound-02}} \\
     & \leq
     \madlocc{\cF}{f}{\frac{d}{2^{J-3}}}{2c}
    a (1 + a^{4-1} + a^{16-1} + \ldots)\mathbbm{1}(J > 1)
    \label{neqn:adaptive-telescoping-bound-04}           \\
     & \leq
     \madlocc{\cF}{f}{\frac{d}{2^{J-3}}}{2c}
    a (1 + a + a^{2} + \ldots)\mathbbm{1}(J > 1)
    \nonumber                                   \\
     & \leq
     \madlocc{\cF}{f}{\frac{d}{2^{J-3}}}{2c}
    \frac{a}{1-a} \mathbbm{1}(J > 1),
    \label{neqn:adaptive-telescoping-bound-05}
\end{align}
where for brevity in \eqref{neqn:adaptive-telescoping-bound-04} we denote
\begin{align*}
    a
    \defined
    \exp\parens{-\frac{n L(\alpha, \beta, C) d^2}{2^{2(J-1)}(C+1)^2}}.
\end{align*}
Since $C$ is assumed to satisfy \eqref{neqn:log-likelihood-bernstein-02a}, and
$L(\alpha, \beta, C)$ is defined as per
\eqref{neqn:log-likelihood-bernstein-02b}, it follows that $a < 1$. Note here
that the above bound \eqref{neqn:adaptive-telescoping-bound-05} holds, provided
that $\PP(\stcomp{A}_{j}) > 0$ for $j < J$ as required by
\Cref{nlem:telescoping-sum-cond-prob}. Suppose that the $\rhs$ of
\eqref{neqn:adaptive-telescoping-bound-05} is strictly smaller than $1$. In that
case for all $j$, $\PP(\stcomp{A}_{j}) > 0$ since bound
\eqref{neqn:adaptive-telescoping-bound-05} holds inductively for all $\PP(A_j)$
for $j \leq J$. On the other hand, if the $\rhs$ of
\eqref{neqn:adaptive-telescoping-bound-05} is $\geq 1$ then
\eqref{neqn:adaptive-telescoping-bound-05} trivially holds. In both cases we
conclude that \eqref{neqn:adaptive-telescoping-bound-05} holds.

If one sets $\varepsilon_J \defined \frac{\sqrt{L(\alpha, \beta, C)}
        d}{2^{(J-1)}(C+1)}$, we have that if
\begin{align*}
    n \varepsilon_J^2 > 2
    \madlocc{\cF}{f}{\varepsilon_J \frac{4(C+1)}{\sqrt{L(\alpha, \beta, C)}}}{2c}
\end{align*} and $a \defined \exp(-n \varepsilon_J^2) <
    1/2 \iff n \varepsilon_J^2 > \log 2$, the above probability in
\eqref{neqn:adaptive-telescoping-bound-05} will be bounded from above by $2
\exp(-n \varepsilon_J^2 / 2)$. This condition is implied when
\begin{align}\label{neqn:adaptive-suff-condition-epsJ}
    n \varepsilon_J^2
    >
    2 \madlocc{\cF}{f}{\varepsilon_J \frac{4(C+1)}{\sqrt{L(\alpha, \beta, C)}}}{2c}
    \vee \log 2.
\end{align}
We now have
\begin{align}\label{neqn:adaptive-mu-upislon-ineq}
    \normb{\nu^*_{\overline{J}} - f}_2
    \leq
    \normb{\Upsilon_{\overline{J}} - \Upsilon_J}_2 +
    \normb{\Upsilon_J - f}_2
    \leq
    3 \varepsilon_J \frac{C+1}{\sqrt{L(\alpha, \beta, C)}},
\end{align} with probability at least $1 - 2\exp(-n\varepsilon_J^2 / 2)$ which holds for all $J$ satisfying \eqref{neqn:adaptive-suff-condition-epsJ}. Here we want to clarify that the last inequality in \eqref{neqn:adaptive-mu-upislon-ineq} follows from the fact that $\|\Upsilon_{\overline{J}} - \Upsilon_J\|_2 \leq d/2^{J-2}$, as per \Cref{nlem:upsilon-is-cauchy-sequence} (and since $\overline{J} \geq J$).
Let $J^*$ be selected as the maximum integer $J$ such that
\eqref{neqn:adaptive-suff-condition-epsJ} holds, or otherwise if such $J$ does
not exist $J^* = 1$. Let $\eta = 3\frac{C+1}{\sqrt{L(\alpha, \beta, C)}}$,
$\underline{C} = 2$ and $C' = 1 / 2$. Observe that by the definition of $J^*$ it
follows that all packing sets encountered prior $J^*$, will have been finite
packing sets. We have established that the following bound holds
\begin{equation*}
    \mathbb{P}(\normb{f - \nu^*_{\overline{J}}}_{2}
    >
    \eta \varepsilon_J)
    \leq
    \underline C \exp(-C'n\varepsilon_J^2) \mathbbm{1}(J > 1)
    \leq
    \underline C \exp(-C'n\varepsilon_J^2) \mathbbm{1}(J^* > 1),
\end{equation*}
for all $1 \leq J \leq J^*$, where this bound also holds in the case when $J^* =
    1$ by exception. Observe that we can extend this bound to all $J \in
    \mathbb{Z}$ and $J \leq J^*$, since for $J < 1$ we have $\eta \varepsilon_J
    \geq 6 d$ and so
\begin{equation*}
    \mathbb{P}(\normb{f - \nu^*_{\overline{J}}}_{2}
    >
    \eta \varepsilon_J)
    \leq
    0
    \leq
    \underline C \exp(-C'n\varepsilon_J^2) \mathbbm{1}(J^* > 1).
\end{equation*}
We conclude that
\begin{equation*}
    \mathbb{P}(\normb{f - \nu^*_{\overline{J}}}_{2}
    >
    \eta \varepsilon_J)
    \leq
    0
    \leq
    \underline C \exp(-C'n\varepsilon_J^2) \mathbbm{1}(J^* > 1),
\end{equation*}
for any $J \leq J^*$. Now for any $\varepsilon_{J-1} > x \geq \varepsilon_{J}$
for $J \leq J^*$ we have that
\begin{align*}
    \mathbb{P}(\normb{f - \nu^*_{\overline{J}}}_{2}
    >
    2\eta x)
     & \leq
    \mathbb{P}(\normb{f - \nu^*_{\overline{J}}}_{2}
    >
    \eta \varepsilon_{J - 1})                                       \\
     & \leq
    \underline C \exp(-C'n\varepsilon_{J-1}^2) \mathbbm{1}(J^* > 1) \\
     & \leq
    \underline C \exp(-C'nx^2)\mathbbm{1}(J^* > 1),
\end{align*}
where the last inequality follows due to the fact that the map $x \mapsto
    \underline C \exp(-C'nx^2)$ is monotonically decreasing for positive reals.
    We will now integrate the tail bound:
\begin{equation}\label{neqn:adaptive-important-prob-bound}
    \mathbb{P}(\normb{f - \nu^*_{\overline{J}}}_{2}
    >
    2 \eta x)
    \leq
    \underline C \exp(-C'nx^2) \mathbbm{1}(J^* > 1),
\end{equation}
which holds true for $x \geq \varepsilon^*$, where $\varepsilon_J =
    \frac{\sqrt{L(\alpha, \beta, C)} d}{2^{(J-1)}(C+1)}$, always (since even if
    $J^* = 1$ by exception, this bound is still valid). We then have
\begin{align*}
    \EE \normb{f - \nu^*_{\overline{J}}}_{2}^{2}
     & =
    \int_{0}^{\infty} 2 x
    \mathbb{P}(\normb{f - \nu^*_{\overline{J}}}_{2} > x) \dlett{x}                                                            \\
     & \leq
    C''' \varepsilon^{*2} + \int_{2\eta\varepsilon^*}^{\infty} 2 x \underline C \exp(-C''nx^2) \mathbbm{1}(J^* > 1) \dlett{x} \\
     & =
    C''' \varepsilon^{*2} + C^{''''}n^{-1}\exp(-C'''''n\varepsilon^{*2})\mathbbm{1}(J^* > 1).
\end{align*}
Now $n\varepsilon^{*2}$ is bigger than a constant (\ie, $\log 2$) otherwise $J^*
    = 1$. Hence, the above is smaller than $\bar C \varepsilon^{*2}$ for some
    absolute constant $\bar C$.
\eprfof

\begin{restatable}[Early stopping in adaptive
estimation]{nrmk}{earlystoppingadaptiveestimator}\label{nrmk:early-stopping-adaptive-estimator}
Suppose that in our adaptive estimation, that we traverse the maximal packing
set tree construction and encounter a density $\Upsilon_{i}$, such that the
cardinality of the set of its children densities is countably infinite, \ie,
$\absa{P_{\Upsilon_{i}}} = \infty$. Then we can simply return
$\nu^{*}(\combdatavec) = \Upsilon_{i}$, in such a case. The reason for this is
that the index $i$ will be necessarily at least equal to $J^{*}$ as defined in
\eqref{neqn:adaptive-upper-bound-suff-cond}, which is what is required for
\eqref{neqn:adaptive-mu-upislon-ineq} to hold.
\end{restatable}

\clearpage
\section{Proofs of
\texorpdfstring{\Cref{sec:examples}}{\autoref{sec:examples}}}\label{appendix:proofs-examples}
\subsection{Formal justification for
\texorpdfstring{\Cref{nexa:lipschitz-density-class}}{\autoref{nexa:lipschitz-density-class}}}\label{app:nexa:lipschitz-density-class}

Before proving
\Cref{nexa:lipschitz-density-class,nexa:bdd-total-variation-density-class,nexa:quad-functional-density-class},
we first prove a useful lemma. This lemma will provide a sufficient condition to
ensure that $L_{2}$-\emph{local} and $L_{2}$-\emph{global} metric entropies are
of the same order for various forms of the density class $\cF$, as specified in
our chosen examples.

\begin{restatable}[Asymptotic order global metric
        entropy]{nlem}{localglobalentropyequivcF}\label{nlem:local-global-ent-cF}
        Let $\cF \subset \cFalphabetaB$, such that for any fixed $\eta > 0$, we
        have $0 < \varepsilon \mapsto \log{\mgloc{\cF}{\varepsilon}} \asymp
        \varepsilon^{- 1 / \eta}$. Then there exists a $c > 0$, such that the
        following holds
    \begin{equation}\label{neqn:local-global-ent-cF-01}
        \log{\mgloc{\cF}{\varepsilon / c}} - \log{\mgloc{\cF}{\varepsilon}}
        \asymp
        \log{\mgloc{\cF}{\varepsilon / c}}
    \end{equation}
\end{restatable}

\bprfof{\Cref{nlem:local-global-ent-cF}} We firstly note that
\eqref{neqn:local-global-ent-cF-01} has the following equivalence
\begin{align}
     & \log{\mgloc{\cF}{\varepsilon / c}} - \log{\mgloc{\cF}{\varepsilon}}
    \asymp
    \log{\mgloc{\cF}{\varepsilon / c}}
    \nonumber                                                              \\
    \iff
    \exists 0 < k_{1} < k_{2} \text{ s.t. }
    k_{1} \log{\mgloc{\cF}{\varepsilon / c}}
    \leq
     & \log{\mgloc{\cF}{\varepsilon / c}} - \log{\mgloc{\cF}{\varepsilon}}
    \leq
    k_{2} \log{\mgloc{\cF}{\varepsilon / c}}
    \label{neqn:local-global-ent-cF-02}
\end{align}
In general, for \Cref{neqn:local-global-ent-cF-02} we observe that since
$\log{\mgloc{\cF}{\varepsilon / c}} > 0$, it follows that
$\log{\mgloc{\cF}{\varepsilon / c}} - \log{\mgloc{\cF}{\varepsilon}} \leq
\log{\mgloc{\cF}{\varepsilon / c}}$. So taking $k_{2} = 1$ will always suffice
to ensure \eqref{neqn:local-global-ent-cF-02} holds. It remains to check that we
can also find a $k_{1} \in (0, 1)$ such that \eqref{neqn:local-global-ent-cF-02}
also holds. In our case, since $\log{\mgloc{\cF}{\varepsilon}} \asymp
\varepsilon^{- 1 / \eta}$ by assumption, we have that
$\log{\mgloc{\cF}{\varepsilon / c}} \geq C_{1} (\varepsilon / c)^{- 1 / \eta}$
and $\log{\mgloc{\cF}{\varepsilon}} \leq C_{2} \varepsilon^{- 1 / \eta}$ for
some universal constants $C_{1}, C_{2} > 0$. It then follows
\begin{align}
    \frac{\log{\mgloc{\cF}{\varepsilon / c}} - \log{\mgloc{\cF}{\varepsilon}}}{\log{\mgloc{\cF}{\varepsilon / c}}}
     & \geq 1 - \parens{\frac{C_{2}}{C_{1}}}c^{- \frac{1}{\eta}}
    \nonumber                                                    \\
     & \geq
    k_{1}
    \tag{as required, for $k_{1} \in (0, 1)$}                    \\
     & \defines
    1 - \delta
    \tag{for some $\delta \in (0, 1)$, since $k_{1} \in (0, 1)$} \\
    \text{if }
    c
     & \geq
    \parens{\frac{C_{2}}{C_{1} \delta}}^{\eta}.
    \label{neqn:local-global-ent-cF-03}
\end{align}
That is, there exists such a $k_{1} \in (0, 1)$, if we choose $c \geq
    \parens{\frac{C_{2}}{C_{1} \delta}}^{\eta}$, for each $\eta > 0$. So indeed
    \eqref{neqn:local-global-ent-cF-01} holds, for the specified class $\cF$, as
    required.
\eprfof

\exalipschitzdensityclass*
\bprfof{\Cref{nexa:lipschitz-density-class}} In order to establish the minimax
rate for $\lipschitzpsi$, we need to show that $\lipschitzpsi$ is a convex
density class, and that there exists a density $f_{\alpha} \in \lipschitzpsi$
that is strictly positively bounded away from 0. We can then apply
\Cref{nprop:extend-zero-bounded-densities}. We first verify that $\lipschitzpsi$
here is a convex density class. To that end, let $f, g \in \lipschitzpsi$, and
let $\kappa \in [0, 1]$, be arbitrary. Then for each $x \in B \defined [0, 1]$,
we observe that
\begin{align}
    (\kappa f + (1 - \kappa) g)(x)
    \defined \kappa f(x) + (1 - \kappa) g(x)
     & \geq \kappa (0) + (1 - \kappa) (0)
    =
    0
    \label{neqn:lipschitz-density-class-convex-01} \\
    (\kappa f + (1 - \kappa) g)(x)
    \defined \kappa f(x) + (1 - \kappa) g(x)
     & \leq
    \kappa \Psi + (1 - \kappa) \Psi
    =
    \Psi
    \label{neqn:lipschitz-density-class-convex-02}
\end{align}
From \eqref{neqn:lipschitz-density-class-convex-01} and
\eqref{neqn:lipschitz-density-class-convex-02}, it follows that
\begin{equation}\label{neqn:lipschitz-density-class-02}
    \kappa f + (1 - \kappa) g \colon B \to [0, \Psi].
\end{equation} Moreover, since $\int_{B} f \dlett{\mu}
    = \int_{B} g \dlett{\mu} = 1$, we have
\begin{equation}\label{neqn:lipschitz-density-class-03}
    \int_{B} (\kappa f + (1 - \kappa) g) \dlett{\mu}
    = \kappa \int_{B} f \dlett{\mu} +
    (1 - \kappa) \int_{B} g \dlett{\mu}
    = 1.
\end{equation}
Since $f, g \in \lipschitzpsi$, we have both $\normb{f}_{q}, \normb{g}_{q} \leq
    \Psi$. Then by the triangle inequality it follows
\begin{equation}\label{neqn:lipschitz-density-class-04}
    \normb{\kappa f + (1 - \kappa) g}_{q}
    \leq
    \normb{\kappa f}_{q} +
    \normb{(1 - \kappa) g}_{q}
    \leq
    \kappa \Psi + (1 - \kappa) \Psi
    =
    \Psi.
\end{equation}
Since $f, g$ are measurable functions, then so is their convex combination, \ie,
$\kappa f + (1 - \kappa) g$. Now we observe
\begin{align}
     & \normb{(\kappa f + (1 - \kappa) g)(x + h)
        - (\kappa f + (1 - \kappa) g)(x)}_{q}
    \nonumber                                    \\
     & =
    \normb{\kappa (f(x + h) - f(x))
        + (1 - \kappa) (g(x + h) - g(x))}_{q}
    \nonumber                                    \\
     & \leq
    \normb{\kappa (f(x + h) - f(x))}_{q} +
    \normb{(1 - \kappa) (g(x + h) - g(x))}_{q}
    \tag{by the triangle inequality.}            \\
     & \leq
    \kappa h^{\gamma} + (1 - \kappa) h^{\gamma}
    \tag{since $f, g \in \lipschitzpsi$}         \\
     & =
    h^{\gamma}
    \label{neqn:lipschitz-density-class-05},
\end{align}
as required. Combining
\eqref{neqn:lipschitz-density-class-02},\eqref{neqn:lipschitz-density-class-03},\eqref{neqn:lipschitz-density-class-04},
and \eqref{neqn:lipschitz-density-class-05} we have shown that $\kappa f + (1 -
\kappa) g \in \lipschitzpsi$. This proves the convexity of $\lipschitzpsi$, as
required.

Now let $f_{\alpha} \sim \distUnif{[B]}$, \ie $f_{\alpha}(x) \defined \Indb{[0,
                1]}{x}$. Therefore, $\normb{f_{\alpha}(x + h) - f(x)}_{q} = 0
                \leq \Psi h^{\gamma}$, for each $x \in B$, and $h > 0$ such that
                $f_{\alpha}(x + h)$ is defined. Moreover,
                $\normb{f_{\alpha}(x)}_{q} = 1 < \Psi$, by assumption, for each
                $1 \leq q \leq \infty$. Now we have that $\int_{B} f_{\alpha}
                \dlett{\mu} = \int_{B} \Indb{[0, 1]}{x} \dlett{\mu}(x) = 1$, and
                $f$ is measurable since it is  a simple function. So indeed we
                have found $f_{\alpha} \in \lipschitzpsi$, such that it is
                $\alpha$-lower bounded (with $\alpha = 1$).

We now proceed to check that $L_{2}$-\emph{global} metric entropy is of the same
order as the $L_{2}$-\emph{local} metric entropy for $\lipschitzpsi$. That is,
we want to check that \eqref{neqn:local-global-ent-cF-01} holds. Here $\cF =
\lipschitzpsi$, with $0 < \varepsilon \mapsto \log{\mgloc{\cF}{\varepsilon}}
\asymp \varepsilon^{- 1 / \gamma}$. Thus, we can apply
\Cref{nlem:local-global-ent-cF} with $\eta \defined \gamma \in (0, 1]$ to
conclude that indeed \eqref{neqn:local-global-ent-cF-01} holds, as required.

Since we have checked all the sufficient conditions in order to apply
\Cref{nprop:extend-zero-bounded-densities} for $\lipschitzpsi$, we can obtain
the minimax rate of density estimation by solving
\begin{equation}
    n \varepsilon^{2}
    \asymp
    \varepsilon^{- \frac{1}{\gamma}}
    \iff
    \varepsilon
    \asymp
    n^{-\frac{\gamma}{2 \gamma + 1}}
    \iff
    \varepsilon^{2}
    \asymp
    n^{-\frac{2 \gamma}{2 \gamma + 1}}.
\end{equation}
So the minimax rate is (up to constants) the order of $n^{-\frac{2 \gamma}{2
                \gamma + 1}}$ as required.
\eprfof

\subsection{Formal justification for
    \texorpdfstring{\Cref{nexa:bdd-total-variation-density-class}}{\autoref{nexa:bdd-total-variation-density-class}}}\label{app:nexa:bdd-total-variation-density-class}

\exabddtotvardensityclass*

\bprfof{\Cref{nexa:bdd-total-variation-density-class}} In order to establish the
minimax rate for $\totbddvzeta$, we need to show that $\totbddvzeta$ is a convex
density class, and that there exists a density $f_{\alpha} \in \totbddvzeta$
that is strictly positively bounded away from 0. We can then apply
\Cref{nprop:extend-zero-bounded-densities}. We first verify that $\totbddvzeta$
here is a convex density class. To that end, let $f, g \in \totbddvzeta$, and
let $\kappa \in [0, 1]$, be arbitrary. Then for each $x \in B \defined [0, 1]$,
it follows by an identical argument to
\eqref{neqn:lipschitz-density-class-convex-01} and
\eqref{neqn:lipschitz-density-class-convex-02} that
\begin{equation}\label{neqn:bdd-total-variation-convex-03}
    \kappa f + (1 - \kappa) g \colon B \to [0, \zeta].
\end{equation} Moreover, since $\int_{B} f \dlett{\mu}
    = \int_{B} g \dlett{\mu} = 1$, we have
\begin{equation}\label{neqn:bdd-total-variation-convex-04}
    \int_{B} (\kappa f + (1 - \kappa) g) \dlett{\mu}
    = \kappa \int_{B} f \dlett{\mu} +
    (1 - \kappa) \int_{B} g \dlett{\mu}
    = 1.
\end{equation}
Since $f, g \in \totbddvzeta$, we have both $\normb{f}_{\infty},
    \normb{g}_{\infty} \leq \zeta$. Then by the triangle inequality it follows
\begin{equation}\label{neqn:bdd-total-variation-convex-05}
    \normb{\kappa f + (1 - \kappa) g}_{\infty}
    \leq
    \normb{\kappa f}_{\infty} +
    \normb{(1 - \kappa) g}_{\infty}
    \leq
    \kappa \zeta + (1 - \kappa) \zeta
    =
    \zeta.
\end{equation}
Since $f, g$ are measurable functions, then so is their convex combination, \ie,
$\kappa f + (1 - \kappa) g$. Finally, fix any $m \in \NN$, and let $a \leq x_{1}
< \cdots < x_{m} \leq b$ be any fixed partition of $B$. Now we observe
\begin{align}
     & \sum_{i=1}^{m - 1} \absb{(\kappa f + (1 - \kappa) g)\parens{x_{i+1}}-(\kappa f + (1 - \kappa) g)\parens{x_{i}}}
    \nonumber                                                                                                          \\
     & =
    \sum_{i=1}^{m - 1} \absb{\kappa (f(x_{i+1}) - f(x_{i}))
        + (1 - \kappa) (g(x_{i+1}) - g(x_{i}))}
    \nonumber                                                                                                          \\
     & \leq
    \kappa \sum_{i=1}^{m - 1} \absb{f(x_{i+1}) - f(x_{i})} +
    (1 - \kappa) \sum_{i=1}^{m - 1} \absb{g(x_{i+1}) - g(x_{i})}
    \tag{by the triangle inequality.}                                                                                  \\
     & \leq
    \kappa V(f) + (1 - \kappa) V(g)
    \nonumber                                                                                                          \\
     & \leq
    \kappa (\zeta) + (1 - \kappa) (\zeta)
    \tag{since $V(f), V(g) \leq \zeta$, by definition of $\totbddvzeta$.}                                              \\
     & =
    \zeta.
    \nonumber
\end{align}
Taking the supremum over all $m \in \NN$ and all partitions of length $m$ of $B$
of the $\lhs$ sum we obtain:
\begin{equation}\label{neqn:bdd-total-variation-convex-06}
    V(\kappa f + (1 - \kappa) g)
    \leq
    \zeta,
\end{equation}
as required. Combining
\eqref{neqn:bdd-total-variation-convex-03},\eqref{neqn:bdd-total-variation-convex-04},\eqref{neqn:bdd-total-variation-convex-05},
and \eqref{neqn:bdd-total-variation-convex-06} we have shown that $\kappa f + (1
- \kappa) g \in \totbddvzeta$. This proves the convexity of $\totbddvzeta$, as
required.

Similar to the proof of \Cref{nexa:lipschitz-density-class}, we let $f_{\alpha}
    \sim \distUnif{[B]}$, \ie $f_{\alpha}(x) \defined \Indb{[0, 1]}{x}$.
    Therefore, $\normb{f}_{\infty} = 1 \leq \zeta$ by assumption. Also, $V(f) =
    0 < \zeta$, by assumption. Now we have that $\int_{B} f_{\alpha} \dlett{\mu}
    = \int_{B} \Indb{[0, 1]}{x} \dlett{\mu}(x) = 1$, and $f$ is measurable since
    it is  a simple function. So indeed we have found $f_{\alpha} \in
    \totbddvzeta$, such that it is $\alpha$-lower bounded (with $\alpha = 1$).

We now proceed to check that $L_{2}$-\emph{global} metric entropy is of the same
order as the $L_{2}$-\emph{local} metric entropy for $\totbddvzeta$. That is, we
want to check that \eqref{neqn:local-global-ent-cF-01} holds. Here $\cF =
\totbddvzeta$, with $0 < \varepsilon \mapsto \log{\mgloc{\cF}{\varepsilon}}
\asymp \varepsilon^{- 1}$. Thus, we can apply \Cref{nlem:local-global-ent-cF}
with $\eta \defined 1$ to conclude that indeed
\eqref{neqn:local-global-ent-cF-01} holds, as required.

Since we have checked all the sufficient conditions in order to apply
\Cref{nprop:extend-zero-bounded-densities} for $\totbddvzeta$, we can obtain the
minimax rate of density estimation by solving
\begin{equation}
    n \varepsilon^{2}
    \asymp
    \varepsilon^{- 1}
    \iff
    \varepsilon
    \asymp
    n^{-\frac{1}{3}}
    \iff
    \varepsilon^{2}
    \asymp
    n^{- \frac{2}{3}}.
\end{equation}
So the minimax rate is (up to constants) the order of $n^{- \frac{2}{3}}$ as
required.
\eprfof

\subsection{Formal justification for
    \texorpdfstring{\Cref{nexa:quad-functional-density-class}}{\autoref{nexa:quad-functional-density-class}}}\label{app:nexa:quad-functional-density-class}
\quadfunctionaldensityclass*

\bprfof{\Cref{nexa:quad-functional-density-class}} In order to establish the
minimax rate for $\quadfuncclass$, we need to show that $\quadfuncclass$ is a
convex density class. We can then apply
\Cref{nprop:extend-zero-bounded-densities}. We first verify that
$\quadfuncclass$ here is a convex density class. To that end, let $f, g \in
\quadfuncclass$, and let $\kappa \in [0, 1]$, be arbitrary. Then for each $x \in
B \defined [0, 1]$, it follows by an identical argument to
\eqref{neqn:lipschitz-density-class-convex-01} and
\eqref{neqn:lipschitz-density-class-convex-02} that
\begin{equation}\label{neqn:quad-functional-density-class-02}
    \kappa f + (1 - \kappa) g \colon B \to [0, \beta].
\end{equation} Moreover, since $\int_{B} f \dlett{\mu}
    = \int_{B} g \dlett{\mu} = 1$, we have
\begin{equation}\label{neqn:quad-functional-density-class-03}
    \int_{B} (\kappa f + (1 - \kappa) g) \dlett{\mu}
    = \kappa \int_{B} f \dlett{\mu} +
    (1 - \kappa) \int_{B} g \dlett{\mu}
    = 1.
\end{equation}
Since $f, g$ are measurable functions, then so is their convex combination, \ie,
$\kappa f + (1 - \kappa) g$. Now we observe
\begin{align}
    \normb{(\kappa f + (1 - \kappa) g)^{\prime \prime}}_{\infty}
    \nonumber
     & =
    \normb{\kappa f^{\prime \prime} + (1 - \kappa) g^{\prime \prime}}_{\infty}
    \tag{by linearity of $2^\textnormal{nd}$ derivative.} \\
     & \leq
    \normb{\kappa f^{\prime \prime}}_{\infty} +
    \normb{(1 - \kappa) g^{\prime \prime}}_{\infty}
    \tag{by the triangle inequality.}                     \\
     & =
    \kappa \normb{f^{\prime \prime}}_{\infty} +
    (1 - \kappa) \normb{g^{\prime \prime}}_{\infty}
    \nonumber                                             \\
     & \leq
    \kappa \gamma + (1 - \kappa) \gamma
    \tag{since $f, g \in \quadfuncclass$}                 \\
     & =
    \gamma
    \label{neqn:quad-functional-density-class-04},
\end{align}
as required. Combining
\eqref{neqn:quad-functional-density-class-02},\eqref{neqn:quad-functional-density-class-03},
and \eqref{neqn:quad-functional-density-class-04} we have shown that $\kappa f +
(1 - \kappa) g \in \quadfuncclass$. This proves the convexity of
$\quadfuncclass$, as required.

Similar to the proof of \Cref{nexa:lipschitz-density-class}, we let $f_{\alpha}
    \sim \distUnif{[B]}$, \ie $f_{\alpha}(x) \defined \Indb{[0, 1]}{x}$. Since
    $\normb{f^{\prime \prime}}_{\infty} = 0 \leq \gamma$. Here, for the boundary
    points of $B \defined [0, 1]$, we are careful to take all derivatives of
    $f_{\alpha}(x)$ at $x = 0$ from the right, and all derivatives from the left
    at $x = 1$. Now we have that $\int_{B} f_{\alpha} \dlett{\mu} = \int_{B}
    \Indb{[0, 1]}{x} \dlett{\mu}(x) = 1$, and $f$ is measurable since it is  a
    simple function. So indeed we have found $f_{\alpha} \in \quadfuncclass$,
    such that it is $\alpha$-lower bounded (with $\alpha = 1$).

We now proceed to check that $L_{2}$-\emph{global} metric entropy is of the same
order as the $L_{2}$-\emph{local} metric entropy for $\quadfuncclass$. That is,
we want to check that \eqref{neqn:local-global-ent-cF-01} holds. Here $\cF =
\quadfuncclass$, with $0 < \varepsilon \mapsto \log{\mgloc{\cF}{\varepsilon}}
\asymp \varepsilon^{- 1 / 4}$. Thus, we can apply
\Cref{nlem:local-global-ent-cF} with $\eta \defined 4$ to conclude that indeed
\eqref{neqn:local-global-ent-cF-01} holds, as required.

Since we have checked all the sufficient conditions in order to apply
\Cref{nprop:extend-zero-bounded-densities} for $\quadfuncclass$, we can obtain
the minimax rate of density estimation by solving
\begin{equation}
    n \varepsilon^{2}
    \asymp
    \varepsilon^{- \frac{1}{2}}
    \iff
    \varepsilon
    \asymp
    n^{-\frac{2}{5}}
    \iff
    \varepsilon^{2}
    \asymp
    n^{- \frac{4}{5}}.
\end{equation}
So the minimax rate is (up to constants) the order of $n^{- \frac{4}{5}}$ as
required.
\eprfof

\subsection{Formal justification for
    \texorpdfstring{\Cref{nexa:convex-mixture-density-class}}{\autoref{nexa:convex-mixture-density-class}}}\label{app:nexa:convex-mixture-density-class}
\exaconvexmixturedensityclass*
\bprfof{\Cref{nexa:convex-mixture-density-class}} Let $\bfG =
    \parens{\bfG_{ij}}_{i, j \in [k]}$ denote the Gram matrix $\bfG_{ij}
    \defined \int_B f_{i} f_{j} \mu(\dlett{x})$. Then it is simple to see that
    for some point $\theta \in \cF$ which can be represented as the convex
    combination $\theta = \sum_{i \in [k]} \alpha_{i} f_i$, the packing set
    should consist of functions $g_{i} = \sum_{j \in [k]} \beta_{ij} f_{j}$
    satisfying both
\begin{align*}
    (\alpha - \beta_i)\T \bfG (\alpha - \beta_i) \leq \varepsilon^2, \\
    (\beta_{i} - \beta_{j})\T \bfG (\beta_{i} - \beta_{j}) > \varepsilon^2/c^2, \mbox{ for } i \neq j,
\end{align*}
where $\beta_i$ are vectors from the $k$-dimensional unit simplex, \ie, $\sum_{j
        \in [k]}\beta_{ij}= 1$, $\beta_{ij} \geq 0$. Now suppose that $\bfG
        \succ 0$. Then upon substituting $\alpha' = \sqrt{\bfG} \alpha$,
        $\beta_i' = \sqrt{\bfG}\beta_{i}$ and dropping the simplex requirements
        on the $\beta$ we obtain the set
\begin{align*}
    \|\alpha ' - \beta_i'\| \leq \varepsilon \\
    \|\beta_i' - \beta_j'\| > \varepsilon/c,
\end{align*}
which is like packing the unit sphere at a distance $1/c$. Hence, the log
cardinality of such a packing is always $\lesssim k$ \cite[see Chapter
5]{wainwright2019high}. If $k$ is not allowed to scale with $n$, we conclude
therefore that the minimax rate is upper bounded by $n^{-1/2}$ which is the
parametric rate as we would expect. If $k$ is allowed to scale with $n$ the rate
is smaller than $\sqrt{\frac{k}{n}}$.
\eprfof

\clearpage
\bibliographystyle{abbrvnat}
\bibliography{refs}

\end{document}